\documentclass[smallextended,envcountsame]{svjour3}

\newif\ifIJCAR

\newif\ifCOLOR
\COLORtrue %

\usepackage{cite}
\usepackage[numbers]{natbib}

\usepackage[utf8]{inputenc}
\usepackage{enumerate}
\usepackage[hang,flushmargin]{footmisc}
\usepackage{amsmath}
\usepackage{amsfonts}
\usepackage{tikz}
\usepackage{booktabs}
\usepackage{xcolor}
\usepackage{bm}
\usepackage{tikz-3dplot}
\usetikzlibrary{calc}
\usepackage{thmtools}
\usepackage{thm-restate}

\declaretheorem[name=Theorem]{thm}
\declaretheorem[name=Proposition]{prop}
\declaretheorem[name=Lemma]{lem}

\ifCOLOR
\usepackage[colorlinks,linkcolor=blue,citecolor=blue,urlcolor=blue]{hyperref}
\newcommand{\ground}{blue}
\else
\usepackage{hyperref}
\newcommand{\ground}{black}
\fi

\definecolor{xred}{rgb}{0.77, 0.12, 0.23}
\definecolor{xgreen}{rgb}{0.3, 0.6, 0}
\definecolor{xblue}{rgb}{0., 0.25, 1}

\newcommand{\ca}{\fill[black] (0,0) circle (3.5pt);}
\newcommand{\cc}{\fill[white] (0,0) circle (3.5pt);}

\ifCOLOR
\newcommand{\cb}{\fill[blue!80!black] (0,0) circle (3.5pt);}
\newcommand{\cd}{\fill[blue!50!white] (0,0) circle (3.5pt);}
\else
\newcommand{\cb}{\fill[black] (-0.11,-0.11) rectangle (0.11,0.11);}
\newcommand{\cd}{\fill[white] (-0.11,-0.11) rectangle (0.11,0.11);}
\fi

\newcommand{\dicea}[6]{
  \draw[fill=white, lightgray] (#1+0.02,#2+0.02) rectangle (#1+.98,#2+.98);
  \begin{scope}[shift={(#1+0.2,#2+0.8)}]{#3}\end{scope}
  \begin{scope}[shift={(#1+0.5,#2+0.8)}]{#4}\end{scope}
  \begin{scope}[shift={(#1+0.8,#2+0.8)}]{#5}\end{scope}
  \begin{scope}[shift={(#1+0.35,#2+0.5)}]{#6}\end{scope}
}

\newcommand{\diceb}[6]{
  \begin{scope}[shift={(#1+0.65,#2+0.5)}]{#3}\end{scope}
  \begin{scope}[shift={(#1+0.2,#2+0.2)}]{#4}\end{scope}
  \begin{scope}[shift={(#1+0.5,#2+0.2)}]{#5}\end{scope}
  \begin{scope}[shift={(#1+0.8,#2+0.2)}]{#6}\end{scope}
}

\newcommand{\drawbox}[4]{
    \pgfmathsetmacro \angle {30}
    \pgfmathsetmacro \xd {{2/3*cos(\angle)}}
    \pgfmathsetmacro \yd {{2/3*sin(\angle)}}
    \pgfmathsetmacro \x {{#1-1+(#2-1)*(\xd)}}
    \pgfmathsetmacro \y {{#3-1+(#2-1)*(\yd)}}

    \draw[fill=#4] (\x,\y) -- (\x+1,\y) -- (\x+1,\y+1) -- (\x,\y+1) -- cycle;
    \draw[fill=#4] (\x,\y+1) -- (\x+\xd,\y+1+\yd) -- (\x+1+\xd,\y+1+\yd) -- (\x+1,\y+1) -- cycle;
    \draw[fill=black!40!#4] (\x+1,\y+1) -- (\x+1+\xd,\y+1+\yd) -- (\x+1+\xd,\y+\yd) -- (\x+1,\y) -- cycle;
}

\newcommand{\drawtop}[4]{
    \pgfmathsetmacro \angle {30}
    \pgfmathsetmacro \xd {{2/3*cos(\angle)}}
    \pgfmathsetmacro \yd {{2/3*sin(\angle)}}
    \pgfmathsetmacro \x {{#1-1+(#2-1)*(\xd)}}
    \pgfmathsetmacro \y {{#3-1+(#2-1)*(\yd)}}

   \draw[fill=#4] (\x,\y+1) -- (\x+\xd,\y+1+\yd) -- (\x+1+\xd,\y+1+\yd) -- (\x+1,\y+1) -- cycle;
}

\newcommand{\drawfront}[4]{
    \pgfmathsetmacro \angle {30}
    \pgfmathsetmacro \xd {{2/3*cos(\angle)}}
    \pgfmathsetmacro \yd {{2/3*sin(\angle)}}
    \pgfmathsetmacro \x {{#1-1+(#2-1)*(\xd)}}
    \pgfmathsetmacro \y {{#3-1+(#2-1)*(\yd)}}

    \draw[fill=#4] (\x,\y) -- (\x+1,\y) -- (\x+1,\y+1) -- (\x,\y+1) -- cycle;
 }

\newcommand{\czero}{c_0}
\newcommand{\cone}{c_1}

\newcommand{\cthree}{c_3}
\newcommand{\la}{\langle}
\newcommand{\ra}{\rangle}

\newcommand{\change}[1]{#1}
\newcommand{\changeprime}[1]{{#1}}

\title{The Resolution of Keller's Conjecture}

\author{Joshua Brakensiek \and Marijn Heule \and John Mackey \and  David Narv\'aez}

\institute{Stanford University, California\and
Carnegie Mellon University, Pennsylvania \and
Carnegie Mellon University, Pennsylvania \and
University of Rochester, New York}

\begin{document}

\maketitle

\begin{abstract}
We consider three graphs, $G_{7,3}$, $G_{7,4}$, and $G_{7,6}$, related to Keller's conjecture in dimension 7. The conjecture is false for this dimension if and only if at least one of the graphs contains a clique of size $2^7 = 128$. We present an automated method to solve this conjecture by encoding the existence of such a clique as a propositional formula. We apply satisfiability solving combined with symmetry-breaking techniques to determine that no such clique exists. This result implies that every unit cube tiling of $\mathbb{R}^7$ contains a facesharing pair of cubes. Since a faceshare-free unit cube tiling of $\mathbb{R}^8$ exists (which we also verify), this completely resolves Keller's conjecture.
\end{abstract}

\section{Introduction}\label{sec:intro}

In 1930, Ott-Heinrich Keller conjectured that any tiling of $n$-dimensional space by translates of the unit cube
must contain a pair of cubes that share a complete ($n-1$)-dimensional face~\cite{keller1930luckenlose}. 
Figure~\ref{fig:dim2} illustrates this for the plane and the $3$-dimensional space.
The conjecture generalized a conjecture of Minkowski~\cite{minkowski} (1907) in which the centers of the cubes were assumed to form a lattice.
Keller's conjecture was proven to be true for $n \leq 6$ by Perron in 1940~\cite{perron1940luckenlose,perron1940luckenlose2}, and in 1942 Haj{\'o}s~\cite{hajos} showed Minkowski's conjecture to be true in all dimensions.

\begin{figure}[b!]
\begin{minipage}{.45\textwidth}
\centering
\begin{tikzpicture}
 \draw (0,0) -- (0,4); \draw (1,0) -- (1,4); \draw (2,0) -- (2,4); \draw (3,0) -- (3,4); \draw (4,0) -- (4,4);
 
\draw[very thick,color=\ground!70!black] (0,0.2) -- (1,0.2) (0,1.2) -- (1,1.2) (0,2.2) -- (1,2.2) (0,3.2) -- (1,3.2);
\draw[very thick,color=\ground!70!black] (1,0.6) -- (2,0.6) (1,1.6) -- (2,1.6) (1,2.6) -- (2,2.6) (1,3.6) -- (2,3.6);
\draw[very thick,color=\ground!70!black] (2,0.4) -- (3,0.4) (2,1.4) -- (3,1.4) (2,2.4) -- (3,2.4) (2,3.4) -- (3,3.4);
\draw[very thick,color=\ground!70!black] (3,0.7) -- (4,0.7) (3,1.7) -- (4,1.7) (3,2.7) -- (4,2.7) (3,3.7) -- (4,3.7);

\end{tikzpicture}
\end{minipage}
\hfill
\begin{minipage}{.45\textwidth}
\centering

\begin{tikzpicture}[scale=1.2]
    \drawbox{0}{1}{1.5}{white}

    \drawbox{1}{2}{1}{white}
    \drawbox{1}{1}{1}{white}

    \drawbox{1}{1.6}{2}{white}

   \drawfront{1}{1.6}{2}{\ground!70!black}

    \drawbox{2}{2.4}{1.6}{white}

    \drawbox{2}{1.4}{1}{white}
    \drawtop{2}{1.4}{1}{\ground!70!black}

    \drawbox{2}{0.4}{0.7}{white}
    \drawtop{2}{0.4}{0.7}{\ground!70!black}

    \end{tikzpicture}
\end{minipage}
~~~
\caption{Left, a tiling of the plane (2-dimensional space) with unit cubes (squares). The bold \ground~edges
are fully face-sharing edges. Right, a partial tiling of the 3-dimensional space with unit cubes. The only way
to tile the entire space would result in a fully face-sharing square at the position of the \ground~squares.}
\label{fig:dim2}
\end{figure}

In 1986 Szab{\'o}~\cite{szabo1986reduction} reduced Keller's conjecture to the study of periodic tilings. 
Using this reduction Corr{\'a}di and Szab{\'o}~\cite{corradiszabo} introduced the Keller graphs: the graph $G_{n,s}$ has vertices $\{0,1,\hdots, 2s-1\}^n$ such that a pair are adjacent if and only if they differ by exactly $s$ in 
at least one coordinate and they differ in at least two coordinates. The size of cliques in $G_{n,s}$ is at most $2^n$~\cite{shor} and the size of the largest clique in $G_{n,s}$ is at most the size of the largest clique in $G_{n,s+1}$.

A clique in $G_{n,s}$ of size $2^n$ demonstrates that Keller's conjecture is false for dimensions greater than or equal to $n$. Lagarias and Shor~\cite{lagarias1992keller} showed that Keller's conjecture is false for $n \geq 10$ in 1992 by exhibiting \change{a} clique of size $2^{10}$ in $G_{10,2}$. In 2002, Mackey~\cite{mackey2002cube} found a clique of size $2^8$ in $G_{8,2}$ to show that Keller's conjecture is false for $n \geq 8$. In 2011, Debroni, Eblen, Langston, Myrvold, Shor, and Weerapurage~\cite{shor} showed that the largest clique in $G_{7,2}$  has size $124$.

In 2015, Kisielewicz and {\L}ysakowska \cite{kisielewicz2017rigid,kisielewicz2015keller}
made substantial progress on reducing the conjecture in dimension $7$. 
More recently, in 2017, Kisielewicz \cite{kisielewicz2017towards} reduced the conjecture in dimension $7$ as follows:
Keller's conjecture is true in dimension $7$ if and only if there does not exist a clique in $G_{7,3}$ of size $2^7$~\cite{lysakowska2019extended}.

The main result of this paper is the following theorem.

\begin{thm}\label{thm:main}
    Neither $G_{7,3}$ nor $G_{7,4}$ nor $G_{7,6}$ contains a clique of size $2^7 = 128$.
\end{thm}

Although proving this property for $G_{7,3}$ suffices to prove Keller's conjecture true in dimension $7$, we also show this for $G_{7,4}$ and $G_{7,6}$ to demonstrate that our methods need only depend on prior work of Kisielewicz and {\L}ysakowska~\cite{kisielewicz2015keller,kisielewicz2017rigid}. In particular, the argument for 
$G_{7,6}$~\cite{kisielewicz2017rigid} predates and is much simpler than the one for $G_{7,4}$~\cite{kisielewicz2015keller}
(although the publication dates indicate otherwise). It is not explicitly stated in either that it suffices to prove that $G_{7,4}$ or $G_{7,6}$ lacks a clique of size $128$ to prove Keller's conjecture.
\ifIJCAR
We show this in the Appendix of the extended version, available at \url{https://arxiv.org/abs/1910.03740}.
\else
We show this in Appendix~\ref{app}.
\fi

We present an approach based on satisfiablity (SAT) solving to show the absence of a clique of size 128. SAT solving has become a powerful tool in computer-aided mathematics in recent years. For example, it was used to prove the Erd\H{o}s discrepancy conjecture with discrepancy 2~\cite{Konev:2015}, the Pythagorean triples problem~\cite{Ptn}, and Schur number five~\cite{S5}. Modern SAT solvers can also emit proofs of unsatisfiability. There exist formally-verified checkers for such proofs as developed in the ACL2, Coq, and Isabelle theorem-proving systems~\cite{Cruz-Filipe2017,Lammich2017}. 

The outline of this paper is as follows. After describing some background concepts in Section~\ref{sec:prelim}, we describe in more detail in Section~\ref{sec:from-conj-to-graph} how Keller's original conjecture on cube tilings reduces to the combinatorial study of Keller graphs. We then present a compact encoding of whether $G_{n,s}$ contains a clique of size $2^n$ as a propositional formula in Section~\ref{sec:encoding}. Without symmetry breaking, these formulas with $n > 5$ are challenging for state-of-the-art tools. However, Keller graphs contain many symmetries. We perform some initial symmetry breaking that is hard to express on the propositional level in Section~\ref{sec:sym}. This allows us to partially fix three vertices. On top of that we add symmetry-breaking clauses in Section~\ref{sec:clsym}. The soundness of their addition has been mechanically verified. 
We prove in Section~\ref{sec:exp} the absence of a clique of size $128$ in $G_{7,3}$, $G_{7,4}$ and $G_{7,6}$. We optimize the proofs of unsatisfiability obtained by the SAT solver and certify them using a formally-verified checker. Finally we draw some conclusions in Section~\ref{sec:conclusions} and present directions for future research. In Appendix~\ref{app}, we fill in some of the technical details \change{from} Section~\ref{sec:from-conj-to-graph}. \change{The Appendix may} serve as the foundation to a fully-formalized proof of Keller's Conjecture.

\begin{figure}[t]
\centering
\begin{tikzpicture}[scale=0.9,
vertex_style/.style={circle,thick, draw},
edge_style/.style={thick, black}]

\node[vertex_style] (0) at (canvas polar cs: radius=0cm,angle=0){$\!\!00\!\!$};
\node[vertex_style] (1) at (canvas polar cs: radius=1.3cm,angle=18+72){$\!\!02\!\!$};
\node[vertex_style] (2) at (canvas polar cs: radius=1.3cm,angle=90+72){$\!\!\textcolor{\ground!90!black}{33}\!\!$};
\node[vertex_style] (3) at (canvas polar cs: radius=1.3cm,angle=162+72){$\!\!\textcolor{\ground!90!black}{3}0\!\!$};
\node[vertex_style] (4) at (canvas polar cs: radius=1.3cm,angle=234+72){$\!\!\textcolor{\ground!90!black}{1}0\!\!$};
\node[vertex_style] (5) at (canvas polar cs: radius=1.3cm,angle=306+72){$\!\!\textcolor{\ground!90!black}{11}\!\!$};

\node[vertex_style] (6) at (canvas polar cs: radius=2.6cm,angle=-18+72){$\!\!2\textcolor{\ground!90!black}{3}\!\!$};
\node[vertex_style] (7) at (canvas polar cs: radius=2.6cm,angle=-90+72){$\!\!\textcolor{\ground!90!black}{3}2\!\!$};
\node[vertex_style] (8) at (canvas polar cs: radius=2.6cm,angle=-162+72){$\!\!22\!\!$};
\node[vertex_style] (9) at (canvas polar cs: radius=2.6cm,angle=-234+72){$\!\!\textcolor{\ground!90!black}{1}2\!\!$};
\node[vertex_style] (10) at (canvas polar cs: radius=2.6cm,angle=-306+72){$\!\!2\textcolor{\ground!90!black}{1}\!\!$};

\node[vertex_style] (11) at (canvas polar cs: radius=3.9cm,angle=18+72){$\!\!20\!\!$};
\node[vertex_style] (12) at (canvas polar cs: radius=3.9cm,angle=90+72){$\!\!0\textcolor{\ground!90!black}{1}\!\!$};
\node[vertex_style] (13) at (canvas polar cs: radius=3.9cm,angle=162+72){$\!\!\textcolor{\ground!90!black}{13}\!\!$};
\node[vertex_style] (14) at (canvas polar cs: radius=3.9cm,angle=234+72){$\!\!\textcolor{\ground!90!black}{31}\!\!$};
\node[vertex_style] (15) at (canvas polar cs: radius=3.9cm,angle=306+72){$\!\!0\textcolor{\ground!90!black}{3}\!\!$};

 \draw[edge_style] (11) -- (7);
 \draw[edge_style] (11) -- (9);
 \draw[edge_style] (12) -- (6);
 \draw[edge_style] (12) -- (8);
 \draw[edge_style] (13) -- (10);
 \draw[edge_style] (13) -- (7);
 \draw[edge_style] (14) -- (6);
 \draw[edge_style] (14) -- (9);
 \draw[edge_style] (15) -- (8);
 \draw[edge_style] (15) -- (10);

 \draw[edge_style] (0) -- (6);
 \draw[edge_style] (0) -- (7);
 \draw[edge_style] (0) -- (8);
 \draw[edge_style] (0) -- (9);
 \draw[edge_style] (0) -- (10);
 
 \draw[edge_style] (1) -- (11);
 \draw[edge_style] (2) -- (12);
 \draw[edge_style] (3) -- (13);
 \draw[edge_style] (4) -- (14);
 \draw[edge_style] (5) -- (15);
 
  \draw[edge_style] (11) -- (12);
  \draw[edge_style] (12) -- (13);
  \draw[edge_style] (13) -- (14);
  \draw[edge_style] (14) -- (15);
  \draw[edge_style] (15) -- (11);
  
 \draw[edge_style] (1) -- (3);
 \draw[edge_style] (3) -- (5);
 \draw[edge_style] (5) -- (2);
 \draw[edge_style] (2) -- (4);
 \draw[edge_style] (4) -- (1);
 
 \draw[edge_style] (1) -- (6);
 \draw[edge_style] (1) -- (10);
 \draw[edge_style] (2) -- (10);
 \draw[edge_style] (2) -- (9);
 \draw[edge_style] (3) -- (9);
 \draw[edge_style] (3) -- (8);
 \draw[edge_style] (4) -- (8);
 \draw[edge_style] (4) -- (7);
 \draw[edge_style] (5) -- (7);
 \draw[edge_style] (5) -- (6);
\end{tikzpicture}
\caption{Illustration of $G_{2,2}$. The coordinates of the vertices are compactly represented by a sequence of the digits.}
\label{g22}
\end{figure}

\section{Preliminaries}\label{sec:prelim}

We present the most important background concepts related to this paper and
introduce some properties of $G_{n,s}$. First, for positive integers $k$, we define two sets: $[k]:=\{1,2,\hdots,k\}$ and $\la k\ra := \{0,1,\hdots , k-1\}$.

\vspace{-2mm}
\paragraph{Keller Graphs.}

The Keller graph $G_{n,s}$ consists of the vertices $\la 2s\ra^n$. Two 
vertices are adjacent if and only if they differ by exactly $s$ in at least 
one coordinate and they differ in at least two coordinates. Figure~\ref{g22} shows a visualization of $G_{2,2}$.

As noted in~\cite{shor}, $\{s w + \la s\ra^n : w \in \{0, 1\}^n \}$ is a partition of the vertices of $G_{n,s}$ into $2^n$ independent sets.
Consequently, any clique in $G_{n,s}$ has at most $2^n$ vertices.
For example, $V(G_{2,2})$ is partitioned as follows:
\begin{eqnarray*}
&&\{2(0,0)+\{0,1\}^2,\,2(0,1)+\{0,1\}^2,\,2(1,0)+\{0,1\}^2,
\,2(1,1)+\{0,1\}^2\} =\\
&&\{\{(0,0), (0,1), (1,0), (1,1)\}, \{(0,2), (0,3), (1,2), (1,3)\},\\
&&\phantom{\{}\{(2,0), (2,1), (3,0), (3,1)\}, \{(2,2), (2,3), (3,2), (3,3)\}\}.
\end{eqnarray*}

We use the above observation to encode whether $G_{n,s}$ has a clique of size $2^n$. Instead of searching for such a clique on the graph representation of $G_{n,s}$, which consists of $(2s)^n$ vertices, we search for $2^n$ vertices, one from each $s w + \la s\ra^n$, such that every pair is adjacent.

For every $i \in \la 2^n\ra$, we let $w(i) = (w_1, w_2,\dots, w_n) \in \{0,1\}^n$ be defined by $i=\sum_{k=1}^n2^{k-1}\cdot w_k$. Given a clique of size $2^n$, we let $c_i$ be its unique element in $s w(i) + \la s\ra^n$ and we let $c_{i,j}$ be the $j$th coordinate of $c_i$.

\paragraph{Useful Automorphisms of Keller Graphs.}

Let $S_n$ be the set of permutations of $[n]$ and let $H_s$ be the set of permutations of $\la 2s\ra$ generated by the swaps $(i\ i+s)$ composed with any permutation of $\la s\ra$ which is identically applied to 
$s + \la s \ra$. \change{For example, when $s = 3$, $(1,2,3,4,5,6)
  \mapsto (3,5,4,6,2,1)$ is such a permutation, where the permutation
  $(a,b,c) \to (c,b,a)$ is applied to both $(1,2,3)$ and $(4,5,6)$ and
  the pairs $(1,4)$ and $(2,5)$ are swapped.} The maps $$(x_1,x_2,\dots,x_n)\mapsto (\tau_1(x_{\sigma(1)}),\tau_2(x_{\sigma(2)}),\dots, \tau_n(x_{\sigma(n)})),$$ where $\sigma \in S_n$ and $\tau_1,\tau_2,\dots,\tau_n \in H_s$ are automorphisms of $G_{n,s}$. Note that applying an automorphism to every vertex of a clique yields another clique of the same size. 

\paragraph{Propositional Formulas.}
We consider formulas in \emph{conjunctive normal form} (CNF), 
which are defined as follows. 
A \emph{literal} is either a variable $x$ (a \emph{positive literal}) 
or the negation $\overline x$ of a variable~$x$ (a \emph{negative literal}). 
The \emph{complement} $\overline l$ of a literal $l$ is defined as 
$\overline l = \overline x$ if $l = x$ and $\overline l = x$ if $l = \overline x$.
For a literal $l$, $\mathit{var}(l)$ denotes the variable of $l$.
A \emph{clause} is a disjunction of literals and a \emph{formula} is a conjunction of clauses.

An \emph{assignment} is a function from a set of variables to the truth values 
1~(\emph{true}) and 0~(\emph{false}).
A literal $l$ is \emph{satisfied} by an assignment $\alpha$ if 
$l$ is positive and \mbox{$\alpha(\mathit{var}(l)) = 1$} or if it is negative and $\alpha(\mathit{var}(l)) = 0$.
A literal is \emph{falsified} by an assignment if its complement is satisfied by the assignment.
A clause is satisfied by an assignment $\alpha$ if it contains a literal that is satisfied by~$\alpha$.
A formula is satisfied by an assignment $\alpha$ if all its clauses are satisfied by $\alpha$.
A formula is \emph{satisfiable} if there exists an assignment that satisfies it and \emph{unsatisfiable} otherwise.

\paragraph{Clausal Proofs.}\label{sec:clausal}

Our proof that Keller's conjecture is true for dimension 7 is predominantly a clausal proof, including a large part of the symmetry breaking. Informally, a clausal proof system allows us to show the unsatisfiability of a CNF formula by continuously deriving more and more clauses until we obtain the empty clause.
Thereby, the addition of a derived clause to the formula and all previously derived clauses must preserve satisfiability. As the empty clause is trivially unsatisfiable, a clausal proof shows the unsatisfiability of the original formula.
Moreover, it must be checkable in polynomial time that each derivation step does preserve satisfiability. This requirement ensures that the correctness of proofs can be efficiently verified. In practice, this is achieved by allowing only the derivation of specific clauses that fulfill some efficiently checkable criterion.

Formally, clausal proof systems are based on the notion of \emph{clause redundancy}. A clause $C$ is \emph{redundant} with respect to a formula $F$ if adding $C$ to $F$ preserves satisfiability. 
Given a formula $F = C_1 \land \dots \land C_m$, a \emph{clausal proof} of $F$
is a sequence $(C_{m+1},\omega_{m+1}), \dots, (C_n, \omega_n)$ of pairs where each $C_i$ is a clause, each $\omega_i$ (called the \emph{witness}) is a string, and $C_n$ is the empty clause~\cite{heule17_shortproofs}. 
Such a sequence gives rise to formulas 
$F_m, F_{m+1}, \dots, F_n$, where $F_i = C_1 \land \dots \land C_i$. 
A clausal proof is \emph{correct} if every clause $C_i$ ($i > m$) is redundant with respect to $F_{i-1}$, and if this redundancy can be checked in polynomial time (with respect to the size of the proof) using the witness $\omega_i$. 

An example for a clausal proof system is the resolution proof system, which only allows the derivation of resolvents (with no or empty witnesses). However, the resolution proof system does not allow to compactly express symmetry breaking.
Instead we will construct a proof in the resolution asymmetric tautology (RAT) proof system. This proof system is also used to validate the results of the SAT competitions~\cite{SC14}. For the details of RAT, we refer to the original paper~\cite{heule17_shortproofs}. Here, we just note that (1)~for RAT clauses, it can be checked efficiently that their addition preserves satisfiability, and (2)~every resolvent is a RAT clause but not vice versa.

\section{From Keller's Conjecture to Keller Graphs}\label{sec:from-conj-to-graph}

In this section, we describe in more detail the connection between Keller graphs and Keller's original conjecture on cube tilings. In particular, we give an overview of the various results since the mid-twentieth century which have contributed toward understanding Keller's conjecture. See Appendix~\ref{app} for a more technical version of this overview (including proofs).

\subsection{Cube Tilings}

Let $d \ge 1$ be the dimension\change{. T}he case $d=7$ is of the \change{primary} interest \change{of} this paper. A \emph{unit cube} (or just \emph{cube}) is a translation of $[0, 1)^d$. In particular, for any $x \in \mathbb R^d$, we define $[0, 1)^d + x$ to be the translated cube $[x_1, x_1+1) \times \cdots [x_d, x_d+1)$. We call $x$ the \emph{corner} of the cube $[0, 1)^d + x$. For brevity, we denote $C^d(x) := [0, 1)^d + x$. We say that two cubes are \emph{disjoint} if they do not intersect. Disjointness is equivalent to the cube corners being ``far apart'' in some coordinate.

\begin{restatable}{prop}{propdisj}
\label{prop:disj}
For all $x, y \in \mathbb R^d$, the cubes $C^d(x)$ and $C^d(y)$ are disjoint if and only if there exists a coordinate $i \in [d]$ such that $|x_i - y_i| \ge 1$.
\end{restatable}

We say that disjoint cubes $C^d(x)$ and $C^d(y)$ are \emph{facesharing} if there exists exactly one coordinate $i \in [d]$ such that $|x_i - y_i| = 1$ and for all $j \in [d]$ such that $j \neq i$, $x_j = y_j$. This is equivalent to saying that $x = y \pm e_i$, where $e_1, \hdots, e_d$ are the standard unit basis vectors.

Let $T \subset \mathbb R^d$ be a set of cube corners. $T$ is a \emph{cube tiling} of $\mathbb{R}^d$ if $[0,1)^d + T = \{C^d(t) : t \in T\}$ is a family of pairwise disjoint cubes such that $\bigcup_{t\in T} C^d(t) = \mathbb{R}^d$. Note that $T = \mathbb Z^d$ produces the standard lattice 
cube tiling of $\mathbb R^d$.

We say that a cube tiling is \emph{faceshare-free} if no pair of distinct cubes in the tiling share a face. For example, $[0, 1)^d + \mathbb Z^d$ is not faceshare-free since, e.g., $[0, 1)^d$ and $[0, 1)^d+e_1$ faceshare. Keller conjectured that all tilings have a facesharing pair of cubes. %

\begin{restatable}[Keller's conjecture~\cite{keller1930luckenlose}]{conj}{kellerconj}
\label{conj:keller}
For all integers $d \ge 1$, there does not exist a faceshare-free 
tiling of $\mathbb R^d$ .
\end{restatable}

\subsection{Structure of Tilings}\label{subsec:structure}

Let $T \subset \mathbb R^d$ be any set of cube corners such that $[0, 1)^d + T$ is a cube tiling. By the definition of cube tiling, we know that for any $x \in \mathbb Z^d$, there exists a unique $t \in T$ such that $x \in C^d(t)$.  The specific values of the coordinates in a tiling are somewhat artificial. That is, the values of the coordinates can be changed while preserving the tiling, as long as the following rule is observed:

\begin{restatable}[``Replacement Lemma'']{lem}{proprepl}
\label{prop:repl}
Let $[0, 1)^d + T$ be a tiling. Fix $a, b \in \mathbb R$ and $i \in [d]$. Define 
\[
T' := \{t \mid t \in T \wedge  t_i \not\equiv a\!\mod 1\} \cup \{t + b e_i \mid t \in T \wedge  t_i \equiv a\!\mod 1\}.
\]
Then, $[0, 1)^d + T'$ is a tiling. Furthermore, if $[0, 1)^d + T$ is faceshare-free and there exists no $t \in T$ such that $t_i \equiv a + b\mod 1$, then $[0, 1)^d + T'$ is faceshare free.
\end{restatable}

\subsection{Reduction to Periodic Tilings}

Next, we show that it suffices to look at \emph{periodic} tilings. We say that $T \subset \mathbb R^d$ and its cube tiling $[0, 1)^d + T$ are \emph{periodic} if for all $t \in T$ and \change{$x \in \mathbb \mathbb Z^d$}, we also have $t + 2x \in T$.  For instance, $T = \mathbb Z^d$ is periodic. In the $t(x)$ notation (see Section~\ref{subsec:structure}), we have that for all $x, y \in \mathbb Z^d$, 
\[
t(x + 2y) = t(x) + 2y.
\]

For a given dimension $d$, if Keller's conjecture is true for all tilings then it is also true for the periodic tilings. The work of Haj\'os~\cite{hajos,hajos1950factorisation} shows that the reverse implication is also true.

\begin{restatable}[Haj\'os~\cite{hajos,hajos1950factorisation}]{thm}{hajos}
\label{thm:hajos}
For all $d \ge 1$, if Keller's conjecture is true for all periodic tilings, then Keller's conjecture is true in dimension $d$.
\end{restatable}

\subsection{Reduction to Keller graphs}

Now we show faceshare-free periodic tilings correspond to cliques in the Keller graph. We say that a periodic tiling is \emph{$s$-discrete} if every coordinate has at most $s$ distinct values modulo $1$ in the tiling. A key observation is that $s$ is bounded by a function of the dimension.

\begin{restatable}[\cite{szabo1986reduction}]{prop}{propszabo}
Every periodic tiling in dimension $d$ is $2^{d-1}$-discrete.
\end{restatable}

As a corollary, we deduce that Keller's conjecture is equivalent to the maxclique problem for a suitable Keller graph.
\begin{restatable}[\cite{corradiszabo}]{thm}{corradiszabo}
\label{thm:corradiszabo}
A faceshare-free periodic tiling exists in dimension $d$ if and only if there exists a clique of size $2^d$ in $G_{d,s}$ where $s = 2^{d-1}$. Therefore, Keller's conjecture in dimension $7$ is equivalent to the lack of a clique of size $128$ in $G_{7,64}$.
\end{restatable}

The relatively recent line of papers by Kisielewicz and {\L}ysakowska~\cite{kisielewicz2017rigid,kisielewicz2015keller,kisielewicz2017towards} improved on the above theorem in dimension $7$ by showing that any potential faceshare-free tilings must be $s$-discrete for $s \in \{3, 4, 6\}$ (in comparison to $64$). Thus, our goal is to rule out the existence of a clique of size $128$ in $G_{7,3}$, $G_{7,4}$, and $G_{7,6}$.

\section[Clique-Existence Encoding]{Clique-Existence Encoding}\label{sec:encoding}

Recall that $G_{n, s}$ has a clique of size $2^n$ if and only if there exist vertices $c_i \in s\change{\cdot} w(i) + \la s \ra^n$ for all $i \in \la 2^n \ra$ such that for all $i \neq i'$ there exist at least two $j \in [n]$ such that $c_{i,j} \neq c_{i',j}$ and there exists at least one $j \in [n]$ such that $c_{i,j} = c_{i',j} \pm s$.

Our CNF will encode the coordinates of the $c_i$. For each $i \in \la 2^n\ra$, $j \in [n]$, $k \in \la s \ra$, we define Boolean variables $x_{i, j, k}$ which are true if and only if $c_{i,j} = s\change{\cdot}w(i)_j + k$. We therefore need to encode that exactly one of $x_{i, j, 0}$, $x_{i, j, 1}$, $\hdots$, $x_{i, j, s-1}$ is true. We use the following clauses
\begin{eqnarray}
\forall i \in \la 2^n\ra, \forall j \in [n], && (x_{i, j, 0} \vee x_{i, j, 1} \vee \cdots \vee x_{i, j, s-1}) \wedge\nonumber%
\bigwedge_{k < k' \in \la s \ra} (\overline{x}_{i, j, k} \vee \overline{x}_{i, j, k'}).\label{eq:cnf1}\\
\end{eqnarray}

Next we enforce that every pair of vertices $c_i$ and $c_{i'}$ in the clique
differ in at least two coordinates. For most pairs of vertices, no clauses are required because $w(i)$ and $w(i')$ differ in at least two positions. Hence, a constraint is only required for two vertices if $w(i)$ and $w(i')$ differ in exactly one coordinate.

Let $\oplus$ be the binary XOR operator and $e_j$ be the indicator vector of the $j$th coordinate.
If $w(i) \oplus w(i') = e_j$, then in order to ensure that
$c_{i}$ and $c_{i'}$ differ in at least two coordinates we need to make sure that there is some coordinate $j' \neq j$ for which $c_{i,j'} \neq c_{i',j'}$
\begin{equation}
    \forall i \neq i' \in \la 2^n\ra \text{ s.t. $w(i) \oplus w(i') = e_j$},\!\!\!\!\!\!\!
    \bigvee_{j' \in [n]\setminus \{j\}, k \in \la s \ra}
    \!\!\!\!\!\!\!(x_{i,j',k} \neq x_{i',j',k}).\label{eq:cnf2}
\end{equation}

We use the Plaisted-Greenbaum~\cite{Plaisted:86} encoding to convert the above constraint into CNF. We refer to the auxiliary variables 
introduced by the encoding as
$y_{i,i',j',k}$, which if true imply $x_{i,j',k} \neq x_{i',j',k}$, or written as an implication
\begin{equation*}
y_{i,i',j',k} \rightarrow (x_{i,j',k} \neq x_{i',j',k})
\end{equation*}
The following two clauses express this implication
\begin{equation}
(\overline y_{i,i',j',k} \lor x_{i,j',k} \lor x_{i',j',k})
\land (\overline y_{i,i',j',k} \lor \overline x_{i,j',k} \lor \overline x_{i',j',k})\label{eq:cnf3}
\end{equation}

Notice that the implication is only in one direction as
Plaisted-Greenbaum takes the polarity of constraints
into account. The clauses that represent the other direction,
i.e., $(y_{i,i',j',k} \lor x_{i,j',k} \lor \overline x_{i',j',k})$ and 
$(y_{i,i',j',k} \lor \overline x_{i,j',k} \lor x_{i',j',k})$ are redundant 
(and more specifically, they are blocked~\cite{Kullmann99}).

Using the auxiliary variables, we can express the constraint (\ref{eq:cnf2})
using clauses of length $s\cdot(n-1)$
\begin{equation}
    \forall i \neq i' \in \la 2^n \ra \text{ s.t. $w(i) \oplus w(i') = e_j$},\!\!\!\!\!\!\!
    \bigvee_{j' \in [n]\setminus \{j\}, k \in \la s \ra}
    \!\!\!\!\!\!\! y_{i,i',j',k}.\label{eq:cnf4}
\end{equation}

The last part of the encoding consists of clauses to ensure that 
each pair of vertices in the clique have at least one coordinate 
in which they differ by exactly $s$.
Observe that 
$c_{i,j} = c_{i',j} \pm s$ implies that $c_{i,j} \neq c_{i',j}$ and  $x_{i,j,k} = x_{i',j,k}$ for all $k \in \la s \ra$.
We use auxiliary variables $z_{i, i', j}$, whose truth implies $c_{i,j} = c_{i',j} \pm s$, or written as
an implication
\begin{eqnarray*}
  &&\!\!\!\!\!\!\!\!\!\!\!\! \forall i \neq i' \in \langle 2^n\rangle, \forall j \in [n]\text{ s.t. }c_{i,j}\neq c_{i',j},\nonumber\\
  &&\!\!\!\!\!\!\!\!
   z_{i, i', j} \rightarrow \big((x_{i,j,0} = x_{i',j,0})
   \land \dots \land (x_{i,j,s-1} = x_{i',j,s-1})\big).\label{eq:cnf5}
\end{eqnarray*}

Notice that the implication is again in one direction only. Below we enforce that some
$z_{i, i', j}$ variables must be true, but there are no constraints that enforce
$z_{i, i', j}$ variables to be false. 

This can be written as a CNF using the following clauses:
\begin{equation}
\bigwedge_{k \in \la s \ra}
\big((\overline z_{i, i', j} \lor {x}_{i,j,k} \lor \overline {x}_{i',j,k}) \land
     (\overline z_{i, i', j} \lor \overline{x}_{i,j,k} \lor {x}_{i',j,k})\big)
\label{eq:cnf6}
\end{equation}

Finally, to make sure that $c_{i,j} = c_{i',j} \pm s$ for some $j \in [n]$, we specify
\begin{equation}
    \forall i \neq i'\in \langle 2^n\rangle,\ \ \ \ \bigvee_{j : c_{i,j} \neq c_{i',j}} {z}_{i, i', j}.\label{eq:cnf7}
\end{equation}

The variables and clauses, including precise formulas for their counts, are summarized in Table~\ref{tbl:cnf-counts}. The sizes of the CNF encodings (before the addition of symmetry breaking clauses) of $G_{7,3}$, $G_{7,4}$, and $G_{7,6}$ are listed in Table~\ref{tbl:cnf-counts2}.  Notice that for fixed $n$, the dependence on $s$ is quadratic, which is better than the $s^{2n}$ dependence one would get in the naive encoding of $G_{n,s}$ as a graph. This compact encoding, when combined with symmetry breaking,  is a core reason why we were able to prove Theorem~\ref{thm:main}.

\begin{table*}[t]
\caption{Summary of variable and clause counts in the CNF encoding.}
\label{tbl:cnf-counts}
\centering
\begin{tabular}{@{}c@{~~~~}c@{~~~~}c}
\toprule
Clauses & New Variable Count & Clause Count\\
\midrule
(\ref{eq:cnf1}) & $2^n \cdot n \cdot s$ & $2^n \cdot n \cdot (1 + \binom{s}{2})$\\
(\ref{eq:cnf3}) & $2^{n-1}\cdot n\cdot s\cdot (n-1)$ & $2^n\cdot n\cdot s\cdot (n-1)$\\
(\ref{eq:cnf4}) &  & $2^{n-1}\cdot n$\\
(\ref{eq:cnf6}) & $2^{2n-2}\cdot n$ & $2^{2n-1} \cdot n \cdot s$\\
(\ref{eq:cnf7}) & & $\binom{2^n}{2}$ \\
\midrule
Total & $2^{n-1}\cdot n\cdot (s(n+1) + 2^{n-1})$ & $2^n\cdot n\cdot\left(\frac{3}{2} + \binom{s}{2} + n\cdot s - s\right) + 2^{2n-1}n s + \binom{2^n}{2}$\\
\bottomrule
\end{tabular}
\end{table*}

\begin{table*}[b]
\caption{Summary of variable and clause counts of the CNF encoding for $G_{7,3}$, $G_{7,4}$, and $G_{7, 6}$. These counts do not include the clauses introduced by the symmetry breaking.}
\label{tbl:cnf-counts2}
\centering
\begin{tabular}{c@{~~~~}c@{~~~~}c}
\toprule
Keller Graph & Variable Count & Clause Count\\
\midrule
$G_{7,3}$ & $39\,424$ & $200\,320$\\
$G_{7,4}$ & $43\,008$ & $265\,728$\\
$G_{7,6}$ & $50\,176$ & $399\,232$\\
\bottomrule
\end{tabular}
\end{table*}
The instances with $n=7$ are too hard for state-of-the-art SAT solvers if no symmetry breaking is applied.
We experimented with general-purpose symmetry-breaking techniques, similar to the symmetry-breaking predicates
produced by shatter~\cite{shatter}. This allows for solving the formula for $G_{7,3}$, but the computation takes a few CPU years. The formulas for $G_{7,4}$ and $G_{7,6}$ with these symmetry-breaking predicates are significantly
harder. 

Instead we employ problem-specific symmetry breaking by making use of the observations in Sections~\ref{sec:sym} and~\ref{sec:clsym}. This allows solving the clique of size $2^n$ existence problem for all three graphs in reasonable time.

\section{Initial Symmetry Breaking}\label{sec:sym}

Our goal is to prove that there exists no clique of size $128$ in $G_{7,s}$ for $s \in \{3, 4, 6\}$. In this section, and the subsequent, we assume that such a clique exists and adapt some of the arguments of Perron~\cite{perron1940luckenlose,perron1940luckenlose2} to show that it may be assumed to have a canonical form.
We will use $\star_i$ to denote an element in $\langle i \rangle$.

\begin{lem}\label{lemma:twocubes}
If there is a clique of size $128$ in $G_{7,s}$, then there is a clique of size $128$ in $G_{7,s}$ containing the vertices $(0, 0, 0, 0, 0, 0, 0)$ and $(s, 1, 0, 0, 0, 0, 0)$.
\end{lem}

\begin{proof} Let $K$ be a clique of size $128$ in $G_{7,s}$. Consider the following sets of vertices in $G_{6,s}$:

$$K_{<s} := \{(v_2,\hdots,v_7) \mid \exists v_1 \in \la s\ra \text { s.t. } (v_1,\hdots,v_7) \in K \}$$ and 
$$K_{\geq s} := \{(v_2,\hdots,v_7) \mid \exists v_1 \in s+\la s\ra \text { s.t. } (v_1,\hdots,v_7) \in K \}.$$

\change{Recall that each vertex of $K$ is of the form $c_i \in s\cdot
  w(i) + \langle s\rangle^7$, $i \in \langle 2^n\rangle.$ Thus,  half
  of the vertices will be in $K_{<s}$ (the ones with $w(i)$ having
  first coordinate $0$) and half will be in $K_{\ge s}$.}
Since $|K_{<s}| + |K_{\geq s}| = 128$, we conclude that $|K_{<s}| = 64$ and $|K_{\geq s}| = 64$.

By the truth of Keller's conjecture in dimension $6$, $K_{<s}$ is not a clique in $G_{6,s}$. Thus, some pair of vertices in $K_{<s}$ are identical in five of the six coordinates. After application of an automorphism, we may without loss of generality assume that this pair is $(s,0,0,0,0,0)$ and $(0,0,0,0,0,0)$. 
Since the pair comes from $K_{<s}$, there exist $v_1 \neq v'_1 \in \la s\ra$ 
such that $(v_1,s,0,0,0,0,0)$ and $(v'_1,0,0,0,0,0,0)$ are in the clique. 

After application of an automorphism that moves $v_1$ to $1$ and $v'_1$ to $0$, we deduce that without loss of generality $(1,s,0,0,0,0,0)$ and $(0,0,0,0,0,0,0)$ are in the clique. Application of the automorphism that interchanges the first two coordinates yields a clique of size $128$ containing the  vertices $\czero = (0, 0, 0, 0, 0, 0, 0)$ and $\cone = (s, 1, 0, 0, 0, 0, 0)$. 
\hfill $\square$
\end{proof}

\smallskip
    
\begin{thm}\label{thm:break}
    If there is a clique of size $128$ in $G_{7,s}$, then there is a clique of size $128$ in $G_{7,s}$ containing the vertices $(0, 0, 0, 0, 0, 0, 0)$, $(s, 1, 0, 0, 0, 0, 0)$, and $(s, s+1, \star_2, \star_2, 1, 1, 1)$.
\end{thm}

\begin{proof}
Using the preceding lemma, we can choose from among all cliques of size $128$ that contain $\czero = (0, 0, 0, 0, 0, 0, 0)$ and 
$\cone = (s, 1, 0, 0, 0, 0, 0)$, one in which $c_3$ has the fewest number of coordinates equal to $0$. Let $\lambda$ be this least number.

Observe that the first two coordinates of $c_3$ must be $(s,s+1)$ in order 
for it to be adjacent with both $c_0$ and $c_1$. Thus, we have 

\[
\begin{array}{ccr@{}c@{\,,\,}c@{\,,\,}c@{\,,\,}c@{\,,\,}c@{\,,\,}c@{\,,\,}c@{}l}
c_0 &=& (&0&0&0&0&0&0&0&)\\
c_1 &=& (&s&1&0&0&0&0&0&)\\
c_3 &=& (&s&s+1&\star_s&\star_s&\star_s&\star_s&\star_s&)
\end{array}
\]

In the above, we can apply automorphisms that fix $0$ in the last five
coordinates to replace $\star_s$ by $\star_2$. We can apply an
automorphism that permutes the last five coordinates to assume that
the $0$'s and $1$'s in $c_3$ are sorted in increasing
order. \change{Imposing such symmetry breaking will be described formally in Section~\ref{sec:clsym}.} Notice that not all of the 
$\star_2$ coordinates in $c_3$ can be $0$, because $\cone$ and $\cthree$ are adjacent and must therefore differ in at least two coordinates.
Hence at least the last coordinate of $c_3$ is $1$.

Case 1) $\lambda = 4$. In this case $c_3 = (s, s+1, 0,0,0,0,1)$. In order for $c_{67}$ to be adjacent with $c_0$, $c_1$, and $c_3$, it must start with $(s,s+1)$ and end with $s+1$:

\[
\begin{array}{ccr@{}c@{\,,\,}c@{\,,\,}c@{\,,\,}c@{\,,\,}c@{\,,\,}c@{\,,\,}c@{}l}
c_0 &=& (&0&0&0&0&0&0&0&)\\
c_1 &=& (&s&1&0&0&0&0&0&)\\
c_3 &=& (&s&s+1&0&0&0&0&1&)\\
c_{67} &=& (&s&s+1&\star_s&\star_s&\star_s&\star_s&s+1&)
\end{array}
\]

Not all $\star_s$ elements in $c_{67}$ can be 0, because $c_3$ and $c_{67}$ differ in at least two coordinates. However, if one of the $\star_s$ elements in $c_{67}$ is nonzero, then we can swap $1$ and $s+1$ in the last coordinate to obtain a clique in which $c_3$ has three or fewer coordinates equal to $0$, contradicting $\lambda = 4$. Thus, $\lambda \leq 3$. 

Case 2) $\lambda = 3$, in which case $c_3 = (s, s+1, 0,0,0,1,1)$:

\[
\begin{array}{ccr@{}c@{\,,\,}c@{\,,\,}c@{\,,\,}c@{\,,\,}c@{\,,\,}c@{\,,\,}c@{}l}
c_0 &=& (&0&0&0&0&0&0&0&)\\
c_1 &=& (&s&1&0&0&0&0&0&)\\
c_3 &=& (&s&s+1&0&0&0&1&1&)\\
c_{35} &=& (&s&s+1&\star_s&\star_s&\star_s&s+1&\star_s&)\\
c_{67} &=& (&s&s+1&\star_s&\star_s&\star_s&\star_s&s+1&)
\end{array}
\]

Since $c_{67}$ is adjacent with $c_0$, $c_1$, and $c_3$, it must start with $(s,s+1)$ and end with $s+1$. Similarly, since $c_{35}$ is adjacent with $c_0$, $c_1$, and $c_3$, it must start with $(s,s+1)$ and have $s+1$ as its penultimate coordinate. Since $c_{35}$ and $c_{67}$ are adjacent, either the last coordinate of $c_{35}$ must be $1$, or the penultimate coordinate of $c_{67}$ must be $1$. Without loss of generality we can assume that the penultimate coordinate of $c_{67}$ is $1$ as we can permute the last two coordinates which would swap $c_{35}$ and $c_{67}$ without involving the other cubes. The remaining three $\star_s$ elements in $c_{67}$ cannot all be $0$, since $c_3$ and $c_{67}$ differ in at least two coordinates. However, if
one of the $\star_s$ elements is non-zero, then we can swap $1$ and $s+1$ in the last coordinate to obtain a clique in which $c_3$ has two or fewer coordinates equal to $0$, contradicting $\lambda = 3$.
Thus, we have $\lambda \leq 2$ and $c_3 = (s, s+1, \star_2,\star_2,1,1,1)$, as desired.
\hfill $\square$
\end{proof}

Notice that most of the symmetry breaking discussed in this section is challenging, if not impossible, to break on the propositional level: The proof of Lemma~\ref{lemma:twocubes} uses the argument that Keller's conjecture holds for dimension $6$, while the proof of Theorem~\ref{thm:break} uses the interchangeability of $1$ and $s+1$, which is not a symmetry on the propositional level. We will break these symmetries by adding some unit clauses to the encoding. All additional symmetry breaking 
will be presented in the next section and will be checked mechanically. 

\section{Clausal Symmetry Breaking}\label{sec:clsym}

Our symmetry-breaking approach starts with enforcing the initial symmetry breaking: We assume that vertices $c_0 = (0,0,0,0,0,0,0)$, $c_1 = (s,1,0,0,0,0,0)$ and $c_3 = (s,s+1,\star_s,\star_s,1,1,1)$ are in our clique $K$, which
follows from Theorem~\ref{thm:break}. We will not use the observation that $\star_s$ occurrences in $c_3$ can be reduced
to $\star_2$ and instead add and validate clauses that realize this reduction.

 We fix the above initial vertices by adding unit clauses to the CNF encoding. This is the only part of the symmetry breaking that is not checked mechanically. Let $\Phi_{7,s}$ be the formula obtained from our encoding in Section~\ref{sec:encoding} together with the unit clauses corresponding to the 19 coordinates fixed among $c_0$, $c_1$ and $c_3$. In this section we will identify several symmetries in $\Phi_{7,s}$ that can be further broken at the CNF level by adding symmetry breaking clauses. \change{We identified these symmetries manually, but they are  present on the CNF level, so we could have used an automated tool such as {\tt shatter}~\cite{shatter}.}
The formula ultimately used in Section~\ref{sec:exp} for the experiments is the result of adding these symmetry breaking clauses to $\Phi_{7,s}$. Symmetry breaking clauses are added in an incremental fashion. For each addition, a clausal proof of its validity with respect to $\Phi_{7,s}$ and the clauses added so far is generated as well. Each of these clausal proofs has been validated using the {\tt drat-trim} proof checker.
 
 Our approach can be described in general terms as identifying groups of coordinates whose assignments exhibit interesting symmetries and calculating the equivalence classes of these assignments. Given a class of symmetric assignments, it holds that one of these assignments can be extended to a clique of size 128 if and only if every assignment in that class can be extended as well. It is then enough to pick a canonical representative for each class, add clauses forbidding every assignment that is not canonical, and finally determine the satisfiability of the formula under the canonical representative of every class of assignments: if no canonical assignment can be extended to a satisfying assignment for the formula, then the formula is unsatisfiable. In order to forbid assignments that are not canonical, we use an approach similar to the one described in~\cite{heule15_symmetry_breaking_drat}.
 
 \subsection{The Last Three Coordinates of $c_{19}$, $c_{35}$ and $c_{67}$}
 
The reasoning in the proof of Theorem~\ref{thm:break} leads to the following forced settings, once we assign $c_3 = (s,s+1,\star_s,\star_s,1,1,1)$ and apply unit propagation:
 
 \begin{itemize}
 \item[] $(c_{19,1},c_{19,2},c_{19,5})=(s,s+1,s+1)$,
 \item[] $(c_{35,1},c_{35,2},c_{35,6})=(s,s+1,s+1)$,
 \item[] $(c_{67,1},c_{67,2},c_{67,7})=(s,s+1,s+1)$.
 \end{itemize}

Let's now focus on the $3\times 3$ matrix of the coordinates below and do a case split on all of the $s^6$ possible assignments of coordinates labeled with~$\star_s$. 

\begin{center}
\begin{tabular}{l|ccc}
&5&6&7\\
\hline
$c_{19}$&$s+1$&$\star_s$&$\star_s$\\
$c_{35}$&$\star_s$&$s+1$&$\star_s$\\
$c_{67}$&$\star_s$&$\star_s$&$s+1$
\end{tabular}
\end{center}
Notice, however, that since the only positions in which $c_{19}$ and $c_{35}$ can differ by exactly $s$ are positions 5 and 6, and since $c_{19,5}$ and $c_{35,6}$ are already set to $s+1$, at least one of $c_{19,6}$ and $c_{35,5}$ has to be set to 1. Similarly, it is not possible for both $c_{35,7}$ and $c_{67,6}$ to not be 1 and for both $c_{67,5}$ and $c_{19,7}$ to not be 1. By the inclusion-exclusion principle, this reasoning alone discards $3(s-1)^2s^4-3(s-1)^4s^2+(s-1)^6$ cases. All of these cases can be blocked by adding the binary clauses
$(x_{19,6,1}\lor x_{35,5,1}) \land 
 (x_{35,7,1}\lor x_{67,6,1}) \land 
 (x_{67,5,1}\lor x_{19,7,1})$.
These three clauses are RAT clauses~\cite{heule12_inprocessing_rules} with respect to the formula $\Phi_{7,s}$.

Furthermore, among the remaining $(2s-1)^3$ cases, several assignment pairs are symmetric. For example, the following two assignments are symmetric because one can be obtained from the other by swapping columns and rows:

\begin{center}
\begin{tabular}{l|ccc}
&5&6&7\\
\hline
$c_{19}$&$s+1$&1&2\\
$c_{35}$&2&$s+1$&2\\
$c_{67}$&1&1&$s+1$
\end{tabular}
\hspace{2cm}
\begin{tabular}{l|ccc}
&5&6&7\\
\hline
$c_{19}$&$s+1$&1&1\\
$c_{35}$&2&$s+1$&1\\
$c_{67}$&2&2&$s+1$
\end{tabular}
\end{center}

As with many problems related to symmetries, we can encode each assignment as a vertex-colored graph and use canonical labeling algorithms to determine a canonical assignment representing all the symmetric assignments of each equivalence class, and which assignments are symmetric to each canonical form. Our approach is similar to the one by McKay and Piperno for isotopy of matrices~\cite{nauty}. %

\begin{figure*}
\centering

\begin{tabular}{c@{~~~~}c@{~~~~}c@{~~~~}c@{~~~~}c}
$(0,1,1,0,0,1)$&
$(0,1,1,0,1,1)$&
$(0,1,1,0,2,1)$&
$(0,1,1,1,0,0)$&
$(0,1,1,1,0,2)$\\
$(0,1,1,1,1,0)$&
$(0,1,1,1,1,1)$&
$(0,1,1,1,1,2)$&
$(0,1,1,1,2,0)$&
$(0,1,1,1,2,1)$\\
$(0,1,1,1,2,2)$&
$(0,1,1,2,1,1)$&
$(0,1,1,2,2,1)$&
$(1,1,0,0,1,1)$&
$(1,1,0,0,2,1)$\\
$(1,1,0,2,1,1)$&
$(1,1,0,2,2,1)$&
$(1,1,1,1,1,1)$&
$(1,1,1,1,1,2)$&
$(1,1,1,1,2,2)$\\
$(1,1,1,2,2,1)$&
$(1,1,2,0,2,1)$&
$(1,1,2,1,2,1)$&
$(1,1,2,1,2,2)$&
$(2,1,1,2,2,1)$\\
$(1,1,2,0,3,1)$&
$(1,1,2,1,3,1)$&
$(1,1,2,1,3,2)$&
\end{tabular}
\caption{The \changeprime{28} canonical cases for \changeprime{$s\ge 4$}.
  \changeprime{The first 25 are the canonical cases for $s=3$} Each vector corresponds to the values of the coordinates $(c_{19,6},c_{19,7},c_{35,5},c_{35,7},c_{67,5},c_{67,6})$.
\change{Note that $c_{19,7} = 1$ in all canonical cases.}}
\label{tbl:s3-new}
\end{figure*}

This additional symmetry breaking reduces the number of cases for the
last three coordinates of the vertices $c_{19}$, $c_{\change{35}}$,
and $c_{67}$ from the trivial $s^6$ to 25 cases for $s=3$ and 28 cases
for $s\geq 4$. Figure~\ref{tbl:s3-new} shows the \changeprime{28}
canonical cases for \changeprime{$s\ge 4$}.

\subsection{Coordinates Three and Four of Vertices $c_3$, $c_{19}$, $c_{35}$ and $c_{67}$}\label{sec:symbr}

\change{From the previous section, we know that either $c_{19,7} = 1$
  or $c_{67,5} = 1$. Via the row/column swap symmetries for the last
  three coordinates of $c_{19},c_{35},$ and $c_{67}$, we can always
  impose without loss of generality that $c_{19,7} = 1$. The particular permutation we use is
      $(c_{19,6},c_{19,7},c_{35,5},c_{35,7},c_{67,5},c_{67,6}) \mapsto
      (c_{67,6},c_{67,5},c_{35,7},c_{35,5},c_{19,7},c_{19,6})$.
   In addition, as long
  as we do not impose any symmetry-breaking conditions on the third
  and fourth coordinates of $c_{19},c_{35},c_{67}$, we can say WLOG
  that} the third and fourth coordinates of $c_{3}$ take values in
$\la 2 \ra$ instead of $\la s \ra$.

\begin{figure}[b]
\centering
$\begin{array}{ccr@{}c@{\,,\,}c@{\,,\,}c@{\,,\,}c@{\,,\,}c@{\,,\,}c@{\,,\,}c@{}l}
c_{0}  &=& (&\bm{ 0}&{\bf 0}&{\bf 0}&{\bf 0}&{\bf 0}&{\bf 0}&{\bf 0}&)\\
c_{1}  &=& (&\bm{ s}&{\bf 1}&{\bf 0}&{\bf 0}&{\bf 0}&{\bf 0}&{\bf 0}&)\\
c_{3}  &=& (&\bm{ s}&\bm{ s+1}&\star_2&\star_2&{\bf 1}&{\bf 1}&{\bf 1}&)\\
c_{19} &=& (&s&s+1&\star_3&\star_3&s+1&\star_3&1&)\\
c_{35} &=& (&s&s+1&\star_4&\star_4&\star_3&s+1&\star_3&)\\
c_{67} &=& (&s&s+1&\star_5&\star_5&\star_4&\star_4&s+1&)
\end{array}$
\caption{Part of the symmetry breaking on the key vertices. The bold coordinates show the (unverified) initial symmetry breaking. The bold $s$ and $s+1$ coordinates in $c_1$ and $c_3$ are also implied by unit propagation.
The additional symmetry breaking is validated by checking a DRAT
proof expressing the symmetry breaking clauses.}
\end{figure}

We break the computation into further cases by enumerating over choices for the third and fourth coordinates of vertices $c_{3}$, $c_{19}$, $c_{\change{35}}$, 
and $c_{67}$.
Up to this point, our description of the partial clique is invariant under the permutations of $\la s-1 \ra$ in the third and fourth coordinates as well as swapping the third and fourth coordinates. With respect to these automorphisms, for $s=3$ there are only $861$ equivalence classes for how to fill in the $\star_s$ cases for these four vertices. For $s=4$ there are $1326$ such equivalence classes, and for $s=6$ there are $1378$ such equivalence classes. This gives a total of $25 \times 861 = 21\,525$ cases to check for $s = 3$, $28 \times 1326 = 37\,128$ cases to check for $s = 4$, and $28 \times 1378 = 38\,584$ cases to check for $s=6$.

\subsection{Identifying Hardest Cases}

In initial experiments we observed for each $s\in\{3,4,6\}$ that out of the 
many thousands of subformulas (cases), one subformula was significantly harder to solve compared to the other subformulas. Figure~\ref{fig:cubes} shows the coordinates of the key vertices of this subformula for $s\in\{3,4,6\}$. Notice that the third and fourth coordinates are all $0$ for all the key vertices. We therefore applied additional symmetry breaking in case all of these coordinates are $0$. Under this case, the third and the fourth coordinates of vertex $c_2$ can be restricted to $(0,0)$, $(0,1)$, and $(1,1)$, and the last three coordinates of $c_2$ can only take values in $\la 3 \ra$. Furthermore, any assignment $(a, b, c)$ to the last three coordinates of $c_2$ is symmetric to the same assignment ``shifted right'', i.e.\ $(c, a, b)$, by swapping columns and rows appropriately. These symmetries define equivalence classes of assignments that can also be broken at the CNF level. Under the case shown in Figure~\ref{fig:cubes}, there are only 33 non-isomorphic assignments remaining
for vertex $c_2$ for $s \geq 3$. We replace the hard case for each $s\in\{3,4,6\}$ by the corresponding 33 cases, thereby
increasing the total number of cases mentioned above by 32. 

\begin{figure}[t]
\hfil
\begin{minipage}{20pt}
$c_{0}$\\
$c_{1}$\\
$c_{3}$\\
$c_{19}$\\
$c_{35}$\\
$c_{67}$
\end{minipage}
\hfil
\begin{minipage}{70pt}
$(0,0,0,0,0,0,0)$\\
$(3,1,0,0,0,0,0)$\\
$(3,4,0,0,1,1,1)$\\
$(3,4,0,0,4,0,1)$\\
$(3,4,0,0,1,4,0)$\\
$(3,4,0,0,0,1,4)$
\end{minipage}
\hfil
\begin{minipage}{70pt}
$(0,0,0,0,0,0,0)$\\
$(4,1,0,0,0,0,0)$\\
$(4,5,0,0,1,1,1)$\\
$(4,5,0,0,5,0,1)$\\
$(4,5,0,0,1,5,0)$\\
$(4,5,0,0,0,1,5)$
\end{minipage}
\hfil
\begin{minipage}{70pt}
$(0,0,0,0,0,0,0)$\\
$(6,1,0,0,0,0,0)$\\
$(6,7,0,0,1,1,1)$\\
$(6,7,0,0,7,0,1)$\\
$(6,7,0,0,1,7,0)$\\
$(6,7,0,0,0,1,7)$
\end{minipage}
\caption{The hardest instance for $s=3$ (left), $s=4$ (middle),
and $s=6$ (right).}
\label{fig:cubes}
\end{figure}

\subsection{SAT Solving}

Each of the cases was solved using a SAT solver, which produced a proof of unsatisfiability that was validated using a formally-verified checker (details are described in the following section). To ensure that the combined cases cover the entire search space, we constructed for each $s\in \{3,4,6\}$ a tautological formula in disjunctive normal form (DNF). The building blocks of a DNF are
conjuctions of literals known as cube. We will use $\alpha$ as a symbol for cubes as they can also be considered variable assignments. For each cube $\alpha$ in the DNF, we prove that the formula after symmetry breaking under $\alpha$ is unsatisfiable. Additionally, we mechanically check that the three DNFs are indeed tautologies.

\begin{figure}[ht]
    \centering
    \includegraphics[width=.9\columnwidth]{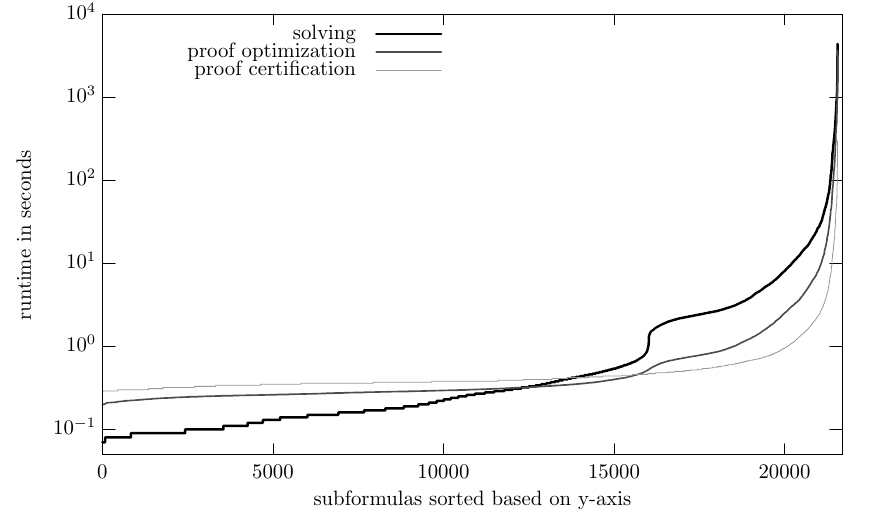}
    \caption{Cactus plot of the runtime in seconds (logscale) to solve the $21\,557$ subformulas of $G_{7,3}$ as well as the times to optimize and certify the proofs of unsatisfiability.}
    \label{fig:cactus3}
\end{figure}

\begin{figure}[ht]
    \centering
    \includegraphics[width=.9\columnwidth]{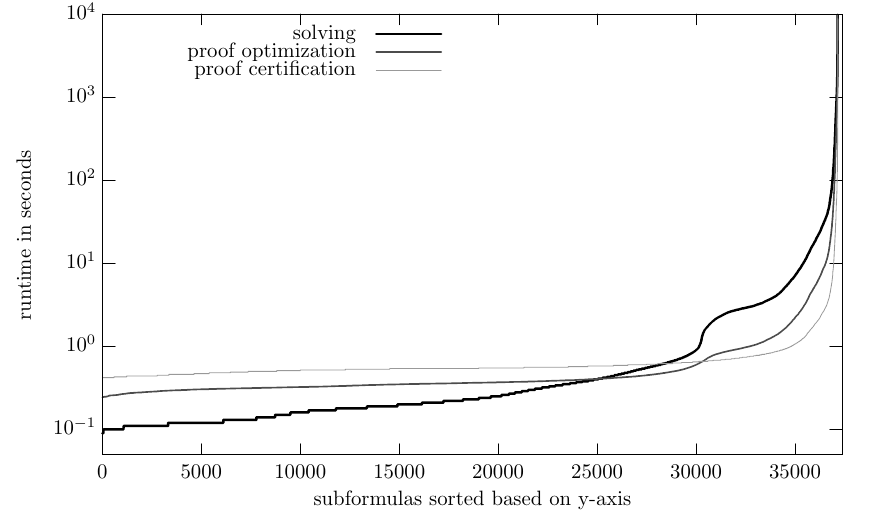}
    \caption{Cactus plot of the runtime in seconds (logscale) to solve the $37\,160$ subformulas of $G_{7,4}$ as well as the times to optimize and certify the proofs of unsatisfiability.}
    \label{fig:cactus4}
\end{figure}

\begin{figure}[ht]
    \centering
    \includegraphics[width=.9\columnwidth]{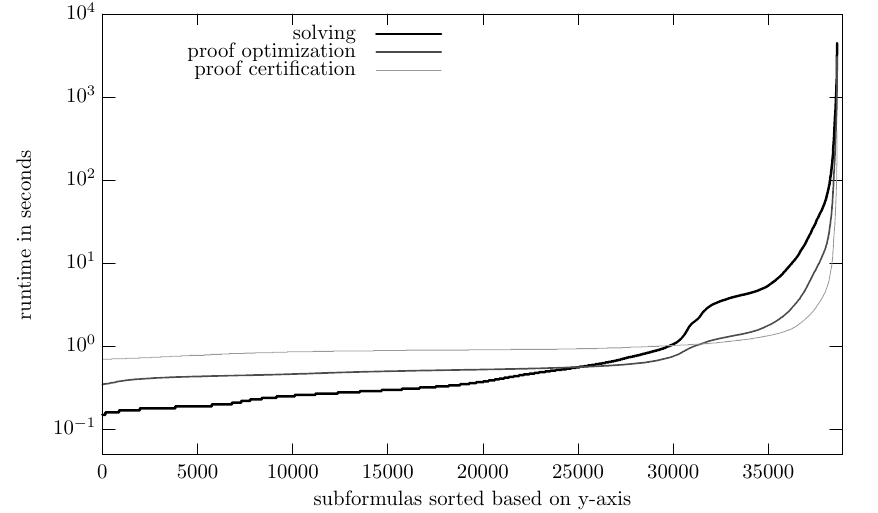}
    \caption{Cactus plot of the runtime in seconds (logscale) to solve the $38\,616$ subformulas of $G_{7,6}$ as well as the times to optimize and certify the proofs of unsatisfiability.}
    \label{fig:cactus6}
\end{figure}

\section{Experiments}\label{sec:exp}

We used the CaDiCaL\footnote{Commit \textsf{92d72896c49b30ad2d50c8e1061ca0681cd23e60} of\\ \url{https://github.com/arminbiere/cadical}} SAT solver developed by  Biere~\cite{cadical} and ran the simulations on a cluster of
Xeon E5-2690 processors with 24 cores per node. CaDiCaL supports proof logging in the DRAT format. We used DRAT-trim~\cite{wetzler2014drat} to optimize the emitted proof of unsatisfiability. Afterwards we certified the optimized proofs with ACL2check, a formally-verified checker~\cite{Cruz-Filipe2017}. 
All of the code that we used is publicly available on GitHub.\footnote{\url{https://github.com/marijnheule/Keller-encode}} We have also made the logs of the computation publicly available on Zenodo.\footnote{\url{https://doi.org/10.5281/zenodo.3755116}}

\begin{table*}[b]
\caption{Summary of solve times for $s = 3, 4, 6$. Times without a unit are in CPU hours. ``No. Hard'' is the number of subformulas which required more than $900$ seconds to solve. ``Hardest'' is the solve time of the hardest subformula in CPU hours,}
\label{tbl:solve-summary}
\centering
\begin{tabular}{@{}c@{~~~~}c@{~~~~}c@{~~~~}c@{~~~~}c@{~~~~}c@{~~~~}c}
\toprule
  $s$ & Tot. Solve & Avg. Solve & Proof Opt. & Proof Cert. & No. Hard & Hardest\\
  \midrule
  3 & 43.27  & 7.23 s & 22.46  & 4.98  & 28 form. & $\approx 1.2$ \\
  4 & 77.00 & 7.46 s & 44.00  & 9.70  & 62 form. & $\approx 2.7$ \\
  6 & 81.85 & 7.63 s & 34.84  & 14.53  & 63 form. & $\approx 1.25$\\
  \bottomrule
\end{tabular}
\end{table*}

\subsection{Results for Dimension 7}

Table~\ref{tbl:solve-summary} summarizes the running times are for experiment. The subformula-solving runtimes for $s = 3, 4$ and $6$ are summarized in cactus plots in Figures~\ref{fig:cactus3}, \ref{fig:cactus4} and \ref{fig:cactus6}.  
The combined size of all unsatisfiability proofs of the subformulas of $s=6$ is
224 gigabytes in the binary DRAT proof format. These proofs contained together $6.18 \cdot 10^9$ proof steps (i.e., additions of redundant clauses). The DRAT-trim proof checker only used $6.39 \cdot 10^8$ proof steps to validate the unsatisfiability of all subformulas. In other words, almost $90\%$ of the clauses generated by CaDiCaL are not required to show unsatisfiability. It is therefore likely that a single DRAT proof for the formula after symmetry breaking can be constructed that is about 20 gigabytes in size. That is significantly smaller compared to other recently solved problems in mathematics that used SAT solvers~\cite{Ptn,S5}.
 
We ran all three experiments simultaneously on 20 nodes on the Lonestar~5 cluster and computing on 24 CPUs per node in parallel.
All instances were reported unsatisfiable and all proofs of unsatisfiability were certified by the formally-verified checker. This proves Theorem~\ref{thm:main}.

\subsection{Refuting Keller's Conjecture in Dimension 8}

To check the accuracy of the CNF encoding, we verified that the generated formulas for $G_{8,2}$, $G_{8,3}$, $G_{8,4}$ and $G_{8,6}$ are satisfiable --- thereby confirming that Keller's conjecture is false for dimension 8. These instances, by themselves, have too many degrees of freedom for the solver to finish. Instead, we added to the CNF the unit clauses consistent with the original clique found in the paper of Mackey~\cite{mackey2002cube} (as suitably embedded for the larger graphs). Specification of the vertices was per the method in Section~\ref{sec:encoding} and~\ref{sec:sym}. These experiments ran on Stanford's Sherlock cluster and took less than a second to confirm satisfiability.

Figure~\ref{fig:clique256} shows an illustration of a clique of size 256 in $G_{8,2}$.
This is the smallest counterexample for Keller's conjecture, both in the dimension ($n=8$)
as in the number of coordinates ($s=2$). 
The illustration uses a 
\ifCOLOR
black, dark blue, white, or light blue dot
\fi
to represent a coordinate set to 0, 1, 2, or 3, respectively.
Notice that for each pair of vertices it holds that they have a complementary 
\ifCOLOR
(black vs white or dark blue vs light blue) dot
\fi
and at least one other different coordinate 
\ifCOLOR
(a different color).
\fi

\begin{figure}[t]

\centering

\begin{tikzpicture}[scale=0.70]

\draw[fill=white,white] (-0.1,-0.1) rectangle (16.1,16.1);

\dicea{0}{0}{\ca}{\ca}{\ca}{\ca}\diceb{0}{0}{\ca}{\ca}{\ca}{\ca}
\dicea{1}{0}{\cc}{\cb}{\ca}{\ca}\diceb{1}{0}{\ca}{\ca}{\ca}{\ca}
\dicea{2}{0}{\ca}{\cc}{\ca}{\cb}\diceb{2}{0}{\cb}{\cb}{\cb}{\ca}
\dicea{3}{0}{\cc}{\cd}{\cb}{\ca}\diceb{3}{0}{\ca}{\cb}{\cb}{\cb}
\dicea{4}{0}{\cb}{\cb}{\cc}{\cb}\diceb{4}{0}{\cb}{\cb}{\ca}{\ca}
\dicea{5}{0}{\cd}{\cb}{\cc}{\cb}\diceb{5}{0}{\cb}{\cb}{\ca}{\cb}
\dicea{6}{0}{\ca}{\cd}{\cc}{\cb}\diceb{6}{0}{\cb}{\cb}{\ca}{\cb}
\dicea{7}{0}{\cc}{\cd}{\cd}{\ca}\diceb{7}{0}{\ca}{\cb}{\ca}{\cb}
\dicea{8}{0}{\cb}{\ca}{\ca}{\cc}\diceb{8}{0}{\cb}{\ca}{\cb}{\cb}
\dicea{9}{0}{\cd}{\ca}{\ca}{\cc}\diceb{9}{0}{\ca}{\ca}{\cb}{\cb}
\dicea{10}{0}{\ca}{\cc}{\cb}{\cd}\diceb{10}{0}{\cb}{\ca}{\cb}{\ca}
\dicea{11}{0}{\cc}{\cc}{\ca}{\cc}\diceb{11}{0}{\cb}{\ca}{\cb}{\cb}
\dicea{12}{0}{\ca}{\ca}{\cc}{\cd}\diceb{12}{0}{\ca}{\ca}{\ca}{\ca}
\dicea{13}{0}{\cc}{\cb}{\cc}{\cd}\diceb{13}{0}{\ca}{\ca}{\ca}{\ca}
\dicea{14}{0}{\ca}{\cc}{\cd}{\cd}\diceb{14}{0}{\cb}{\cb}{\cb}{\ca}
\dicea{15}{0}{\cc}{\cd}{\cc}{\cc}\diceb{15}{0}{\ca}{\cb}{\cb}{\cb}
\dicea{0}{1}{\cb}{\ca}{\ca}{\ca}\diceb{0}{1}{\cc}{\ca}{\ca}{\ca}
\dicea{1}{1}{\cd}{\ca}{\ca}{\cb}\diceb{1}{1}{\cc}{\cb}{\cb}{\cb}
\dicea{2}{1}{\ca}{\cc}{\cb}{\ca}\diceb{2}{1}{\cd}{\cb}{\cb}{\ca}
\dicea{3}{1}{\cc}{\cc}{\ca}{\ca}\diceb{3}{1}{\cc}{\ca}{\ca}{\ca}
\dicea{4}{1}{\ca}{\ca}{\cc}{\cb}\diceb{4}{1}{\cd}{\cb}{\ca}{\cb}
\dicea{5}{1}{\cc}{\cb}{\cc}{\cb}\diceb{5}{1}{\cd}{\cb}{\ca}{\cb}
\dicea{6}{1}{\ca}{\cc}{\cd}{\ca}\diceb{6}{1}{\cd}{\cb}{\ca}{\ca}
\dicea{7}{1}{\cc}{\cd}{\cc}{\cb}\diceb{7}{1}{\cc}{\cb}{\ca}{\cb}
\dicea{8}{1}{\cb}{\cb}{\ca}{\cc}\diceb{8}{1}{\cd}{\ca}{\cb}{\ca}
\dicea{9}{1}{\cd}{\ca}{\cb}{\cd}\diceb{9}{1}{\cc}{\ca}{\cb}{\cb}
\dicea{10}{1}{\cb}{\cd}{\ca}{\cc}\diceb{10}{1}{\cd}{\ca}{\cb}{\cb}
\dicea{11}{1}{\cd}{\cc}{\ca}{\cc}\diceb{11}{1}{\cd}{\ca}{\cb}{\cb}
\dicea{12}{1}{\cb}{\ca}{\cc}{\cd}\diceb{12}{1}{\cc}{\ca}{\ca}{\ca}
\dicea{13}{1}{\cd}{\ca}{\cd}{\cd}\diceb{13}{1}{\cc}{\cb}{\cb}{\cb}
\dicea{14}{1}{\ca}{\cc}{\cc}{\cc}\diceb{14}{1}{\cd}{\cb}{\cb}{\ca}
\dicea{15}{1}{\cc}{\cc}{\cc}{\cd}\diceb{15}{1}{\cc}{\ca}{\ca}{\ca}
\dicea{0}{2}{\ca}{\ca}{\ca}{\cb}\diceb{0}{2}{\ca}{\cc}{\ca}{\ca}
\dicea{1}{2}{\cc}{\cb}{\ca}{\cb}\diceb{1}{2}{\ca}{\cc}{\ca}{\ca}
\dicea{2}{2}{\ca}{\cc}{\cb}{\cb}\diceb{2}{2}{\cb}{\cd}{\cb}{\ca}
\dicea{3}{2}{\cc}{\cd}{\ca}{\ca}\diceb{3}{2}{\ca}{\cd}{\cb}{\cb}
\dicea{4}{2}{\cb}{\ca}{\cc}{\cb}\diceb{4}{2}{\cb}{\cd}{\ca}{\cb}
\dicea{5}{2}{\cd}{\ca}{\cc}{\cb}\diceb{5}{2}{\ca}{\cd}{\ca}{\cb}
\dicea{6}{2}{\ca}{\cc}{\cd}{\ca}\diceb{6}{2}{\cb}{\cd}{\ca}{\ca}
\dicea{7}{2}{\cc}{\cc}{\cc}{\cb}\diceb{7}{2}{\cb}{\cd}{\ca}{\cb}
\dicea{8}{2}{\cb}{\cb}{\cb}{\cd}\diceb{8}{2}{\cb}{\cc}{\cb}{\ca}
\dicea{9}{2}{\cd}{\cb}{\cb}{\cd}\diceb{9}{2}{\cb}{\cc}{\cb}{\cb}
\dicea{10}{2}{\ca}{\cd}{\cb}{\cd}\diceb{10}{2}{\cb}{\cc}{\cb}{\cb}
\dicea{11}{2}{\cc}{\cd}{\ca}{\cc}\diceb{11}{2}{\ca}{\cc}{\cb}{\cb}
\dicea{12}{2}{\ca}{\ca}{\cd}{\cd}\diceb{12}{2}{\ca}{\cc}{\ca}{\ca}
\dicea{13}{2}{\cc}{\cb}{\cd}{\cd}\diceb{13}{2}{\ca}{\cc}{\ca}{\ca}
\dicea{14}{2}{\ca}{\cc}{\cd}{\cc}\diceb{14}{2}{\cb}{\cd}{\cb}{\ca}
\dicea{15}{2}{\cc}{\cd}{\cc}{\cd}\diceb{15}{2}{\ca}{\cd}{\cb}{\cb}
\dicea{0}{3}{\cb}{\ca}{\ca}{\cb}\diceb{0}{3}{\cc}{\cc}{\ca}{\ca}
\dicea{1}{3}{\cd}{\ca}{\cb}{\cb}\diceb{1}{3}{\cc}{\cd}{\cb}{\cb}
\dicea{2}{3}{\ca}{\cc}{\ca}{\ca}\diceb{2}{3}{\cd}{\cd}{\cb}{\ca}
\dicea{3}{3}{\cc}{\cc}{\ca}{\cb}\diceb{3}{3}{\cc}{\cc}{\ca}{\ca}
\dicea{4}{3}{\cb}{\cb}{\cc}{\cb}\diceb{4}{3}{\cd}{\cd}{\ca}{\ca}
\dicea{5}{3}{\cd}{\ca}{\cd}{\ca}\diceb{5}{3}{\cc}{\cd}{\ca}{\cb}
\dicea{6}{3}{\cb}{\cd}{\cc}{\cb}\diceb{6}{3}{\cd}{\cd}{\ca}{\cb}
\dicea{7}{3}{\cd}{\cc}{\cc}{\cb}\diceb{7}{3}{\cd}{\cd}{\ca}{\cb}
\dicea{8}{3}{\ca}{\ca}{\cb}{\cd}\diceb{8}{3}{\cd}{\cc}{\cb}{\cb}
\dicea{9}{3}{\cc}{\cb}{\cb}{\cd}\diceb{9}{3}{\cd}{\cc}{\cb}{\cb}
\dicea{10}{3}{\ca}{\cc}{\ca}{\cc}\diceb{10}{3}{\cd}{\cc}{\cb}{\ca}
\dicea{11}{3}{\cc}{\cd}{\cb}{\cd}\diceb{11}{3}{\cc}{\cc}{\cb}{\cb}
\dicea{12}{3}{\cb}{\ca}{\cd}{\cd}\diceb{12}{3}{\cc}{\cc}{\ca}{\ca}
\dicea{13}{3}{\cd}{\ca}{\cd}{\cc}\diceb{13}{3}{\cc}{\cd}{\cb}{\cb}
\dicea{14}{3}{\ca}{\cc}{\cc}{\cd}\diceb{14}{3}{\cd}{\cd}{\cb}{\ca}
\dicea{15}{3}{\cc}{\cc}{\cd}{\cd}\diceb{15}{3}{\cc}{\cc}{\ca}{\ca}
\dicea{0}{4}{\ca}{\ca}{\cb}{\ca}\diceb{0}{4}{\ca}{\ca}{\cc}{\ca}
\dicea{1}{4}{\cc}{\cb}{\cb}{\ca}\diceb{1}{4}{\ca}{\ca}{\cc}{\ca}
\dicea{2}{4}{\ca}{\cc}{\ca}{\ca}\diceb{2}{4}{\cb}{\cb}{\cd}{\ca}
\dicea{3}{4}{\cc}{\cd}{\cb}{\cb}\diceb{3}{4}{\ca}{\cb}{\cd}{\cb}
\dicea{4}{4}{\cb}{\ca}{\cd}{\ca}\diceb{4}{4}{\cb}{\cb}{\cc}{\cb}
\dicea{5}{4}{\cd}{\ca}{\cd}{\ca}\diceb{5}{4}{\ca}{\cb}{\cc}{\cb}
\dicea{6}{4}{\ca}{\cc}{\cc}{\cb}\diceb{6}{4}{\cb}{\cb}{\cc}{\ca}
\dicea{7}{4}{\cc}{\cc}{\cd}{\ca}\diceb{7}{4}{\cb}{\cb}{\cc}{\cb}
\dicea{8}{4}{\cb}{\cb}{\ca}{\cc}\diceb{8}{4}{\cb}{\ca}{\cd}{\ca}
\dicea{9}{4}{\cd}{\cb}{\ca}{\cc}\diceb{9}{4}{\cb}{\ca}{\cd}{\cb}
\dicea{10}{4}{\ca}{\cd}{\ca}{\cc}\diceb{10}{4}{\cb}{\ca}{\cd}{\cb}
\dicea{11}{4}{\cc}{\cd}{\cb}{\cd}\diceb{11}{4}{\ca}{\ca}{\cd}{\cb}
\dicea{12}{4}{\ca}{\ca}{\cc}{\cc}\diceb{12}{4}{\ca}{\ca}{\cc}{\ca}
\dicea{13}{4}{\cc}{\cb}{\cc}{\cc}\diceb{13}{4}{\ca}{\ca}{\cc}{\ca}
\dicea{14}{4}{\ca}{\cc}{\cc}{\cd}\diceb{14}{4}{\cb}{\cb}{\cd}{\ca}
\dicea{15}{4}{\cc}{\cd}{\cd}{\cc}\diceb{15}{4}{\ca}{\cb}{\cd}{\cb}
\dicea{0}{5}{\cb}{\ca}{\cb}{\ca}\diceb{0}{5}{\cc}{\ca}{\cc}{\ca}
\dicea{1}{5}{\cd}{\ca}{\ca}{\ca}\diceb{1}{5}{\cc}{\cb}{\cd}{\cb}
\dicea{2}{5}{\ca}{\cc}{\cb}{\cb}\diceb{2}{5}{\cd}{\cb}{\cd}{\ca}
\dicea{3}{5}{\cc}{\cc}{\cb}{\ca}\diceb{3}{5}{\cc}{\ca}{\cc}{\ca}
\dicea{4}{5}{\cb}{\cb}{\cd}{\ca}\diceb{4}{5}{\cd}{\cb}{\cc}{\ca}
\dicea{5}{5}{\cd}{\ca}{\cc}{\cb}\diceb{5}{5}{\cc}{\cb}{\cc}{\cb}
\dicea{6}{5}{\cb}{\cd}{\cd}{\ca}\diceb{6}{5}{\cd}{\cb}{\cc}{\cb}
\dicea{7}{5}{\cd}{\cc}{\cd}{\ca}\diceb{7}{5}{\cd}{\cb}{\cc}{\cb}
\dicea{8}{5}{\ca}{\ca}{\ca}{\cc}\diceb{8}{5}{\cd}{\ca}{\cd}{\cb}
\dicea{9}{5}{\cc}{\cb}{\ca}{\cc}\diceb{9}{5}{\cd}{\ca}{\cd}{\cb}
\dicea{10}{5}{\ca}{\cc}{\cb}{\cd}\diceb{10}{5}{\cd}{\ca}{\cd}{\ca}
\dicea{11}{5}{\cc}{\cd}{\ca}{\cc}\diceb{11}{5}{\cc}{\ca}{\cd}{\cb}
\dicea{12}{5}{\cb}{\ca}{\cc}{\cc}\diceb{12}{5}{\cc}{\ca}{\cc}{\ca}
\dicea{13}{5}{\cd}{\ca}{\cc}{\cd}\diceb{13}{5}{\cc}{\cb}{\cd}{\cb}
\dicea{14}{5}{\ca}{\cc}{\cd}{\cc}\diceb{14}{5}{\cd}{\cb}{\cd}{\ca}
\dicea{15}{5}{\cc}{\cc}{\cc}{\cc}\diceb{15}{5}{\cc}{\ca}{\cc}{\ca}
\dicea{0}{6}{\ca}{\ca}{\cb}{\cb}\diceb{0}{6}{\ca}{\cc}{\cc}{\ca}
\dicea{1}{6}{\cc}{\cb}{\cb}{\cb}\diceb{1}{6}{\ca}{\cc}{\cc}{\ca}
\dicea{2}{6}{\ca}{\cc}{\cb}{\ca}\diceb{2}{6}{\cb}{\cd}{\cd}{\ca}
\dicea{3}{6}{\cc}{\cd}{\ca}{\cb}\diceb{3}{6}{\ca}{\cd}{\cd}{\cb}
\dicea{4}{6}{\cb}{\cb}{\cd}{\ca}\diceb{4}{6}{\cb}{\cd}{\cc}{\ca}
\dicea{5}{6}{\cd}{\cb}{\cd}{\ca}\diceb{5}{6}{\cb}{\cd}{\cc}{\cb}
\dicea{6}{6}{\ca}{\cd}{\cd}{\ca}\diceb{6}{6}{\cb}{\cd}{\cc}{\cb}
\dicea{7}{6}{\cc}{\cd}{\cc}{\cb}\diceb{7}{6}{\ca}{\cd}{\cc}{\cb}
\dicea{8}{6}{\cb}{\ca}{\cb}{\cd}\diceb{8}{6}{\cb}{\cc}{\cd}{\cb}
\dicea{9}{6}{\cd}{\ca}{\cb}{\cd}\diceb{9}{6}{\ca}{\cc}{\cd}{\cb}
\dicea{10}{6}{\ca}{\cc}{\ca}{\cc}\diceb{10}{6}{\cb}{\cc}{\cd}{\ca}
\dicea{11}{6}{\cc}{\cc}{\cb}{\cd}\diceb{11}{6}{\cb}{\cc}{\cd}{\cb}
\dicea{12}{6}{\ca}{\ca}{\cd}{\cc}\diceb{12}{6}{\ca}{\cc}{\cc}{\ca}
\dicea{13}{6}{\cc}{\cb}{\cd}{\cc}\diceb{13}{6}{\ca}{\cc}{\cc}{\ca}
\dicea{14}{6}{\ca}{\cc}{\cc}{\cc}\diceb{14}{6}{\cb}{\cd}{\cd}{\ca}
\dicea{15}{6}{\cc}{\cd}{\cd}{\cd}\diceb{15}{6}{\ca}{\cd}{\cd}{\cb}
\dicea{0}{7}{\cb}{\ca}{\cb}{\cb}\diceb{0}{7}{\cc}{\cc}{\cc}{\ca}
\dicea{1}{7}{\cd}{\ca}{\cb}{\ca}\diceb{1}{7}{\cc}{\cd}{\cd}{\cb}
\dicea{2}{7}{\ca}{\cc}{\ca}{\cb}\diceb{2}{7}{\cd}{\cd}{\cd}{\ca}
\dicea{3}{7}{\cc}{\cc}{\cb}{\cb}\diceb{3}{7}{\cc}{\cc}{\cc}{\ca}
\dicea{4}{7}{\ca}{\ca}{\cd}{\ca}\diceb{4}{7}{\cd}{\cd}{\cc}{\cb}
\dicea{5}{7}{\cc}{\cb}{\cd}{\ca}\diceb{5}{7}{\cd}{\cd}{\cc}{\cb}
\dicea{6}{7}{\ca}{\cc}{\cc}{\cb}\diceb{6}{7}{\cd}{\cd}{\cc}{\ca}
\dicea{7}{7}{\cc}{\cd}{\cd}{\ca}\diceb{7}{7}{\cc}{\cd}{\cc}{\cb}
\dicea{8}{7}{\cb}{\cb}{\cb}{\cd}\diceb{8}{7}{\cd}{\cc}{\cd}{\ca}
\dicea{9}{7}{\cd}{\ca}{\ca}{\cc}\diceb{9}{7}{\cc}{\cc}{\cd}{\cb}
\dicea{10}{7}{\cb}{\cd}{\cb}{\cd}\diceb{10}{7}{\cd}{\cc}{\cd}{\cb}
\dicea{11}{7}{\cd}{\cc}{\cb}{\cd}\diceb{11}{7}{\cd}{\cc}{\cd}{\cb}
\dicea{12}{7}{\cb}{\ca}{\cd}{\cc}\diceb{12}{7}{\cc}{\cc}{\cc}{\ca}
\dicea{13}{7}{\cd}{\ca}{\cc}{\cc}\diceb{13}{7}{\cc}{\cd}{\cd}{\cb}
\dicea{14}{7}{\ca}{\cc}{\cd}{\cd}\diceb{14}{7}{\cd}{\cd}{\cd}{\ca}
\dicea{15}{7}{\cc}{\cc}{\cd}{\cc}\diceb{15}{7}{\cc}{\cc}{\cc}{\ca}
\dicea{0}{8}{\cb}{\cb}{\cb}{\ca}\diceb{0}{8}{\cb}{\cb}{\cb}{\cc}
\dicea{1}{8}{\cd}{\cb}{\ca}{\ca}\diceb{1}{8}{\ca}{\ca}{\ca}{\cc}
\dicea{2}{8}{\ca}{\cd}{\ca}{\ca}\diceb{2}{8}{\ca}{\ca}{\ca}{\cc}
\dicea{3}{8}{\cc}{\cd}{\ca}{\cb}\diceb{3}{8}{\ca}{\cb}{\cb}{\cd}
\dicea{4}{8}{\cb}{\cb}{\cd}{\ca}\diceb{4}{8}{\cb}{\cb}{\ca}{\cc}
\dicea{5}{8}{\cd}{\ca}{\cc}{\cb}\diceb{5}{8}{\ca}{\cb}{\ca}{\cd}
\dicea{6}{8}{\cb}{\cd}{\cc}{\cb}\diceb{6}{8}{\cb}{\cb}{\ca}{\cd}
\dicea{7}{8}{\cd}{\cc}{\cc}{\cb}\diceb{7}{8}{\cb}{\cb}{\ca}{\cd}
\dicea{8}{8}{\ca}{\ca}{\ca}{\cc}\diceb{8}{8}{\cb}{\ca}{\cb}{\cd}
\dicea{9}{8}{\cc}{\cb}{\ca}{\cc}\diceb{9}{8}{\cb}{\ca}{\cb}{\cd}
\dicea{10}{8}{\ca}{\cc}{\ca}{\cc}\diceb{10}{8}{\cb}{\ca}{\cb}{\cc}
\dicea{11}{8}{\cc}{\cd}{\cb}{\cd}\diceb{11}{8}{\ca}{\ca}{\cb}{\cd}
\dicea{12}{8}{\cb}{\cb}{\cc}{\cc}\diceb{12}{8}{\cb}{\cb}{\cb}{\cc}
\dicea{13}{8}{\cd}{\cb}{\cc}{\cd}\diceb{13}{8}{\ca}{\ca}{\ca}{\cc}
\dicea{14}{8}{\ca}{\cd}{\cc}{\cd}\diceb{14}{8}{\ca}{\ca}{\ca}{\cc}
\dicea{15}{8}{\cc}{\cd}{\cd}{\cd}\diceb{15}{8}{\ca}{\cb}{\cb}{\cd}
\dicea{0}{9}{\cb}{\cb}{\ca}{\cb}\diceb{0}{9}{\cd}{\cb}{\cb}{\cc}
\dicea{1}{9}{\cd}{\ca}{\cb}{\ca}\diceb{1}{9}{\cc}{\cb}{\cb}{\cd}
\dicea{2}{9}{\cb}{\cd}{\ca}{\ca}\diceb{2}{9}{\cc}{\ca}{\ca}{\cc}
\dicea{3}{9}{\cd}{\cc}{\ca}{\ca}\diceb{3}{9}{\cc}{\ca}{\ca}{\cc}
\dicea{4}{9}{\cb}{\ca}{\cc}{\cb}\diceb{4}{9}{\cd}{\cb}{\ca}{\cd}
\dicea{5}{9}{\cd}{\ca}{\cd}{\ca}\diceb{5}{9}{\cc}{\cb}{\ca}{\cd}
\dicea{6}{9}{\ca}{\cc}{\cc}{\cb}\diceb{6}{9}{\cd}{\cb}{\ca}{\cc}
\dicea{7}{9}{\cc}{\cc}{\cc}{\cb}\diceb{7}{9}{\cd}{\cb}{\ca}{\cd}
\dicea{8}{9}{\cb}{\cb}{\cb}{\cd}\diceb{8}{9}{\cd}{\ca}{\cb}{\cc}
\dicea{9}{9}{\cd}{\cb}{\ca}{\cc}\diceb{9}{9}{\cd}{\ca}{\cb}{\cd}
\dicea{10}{9}{\ca}{\cd}{\ca}{\cc}\diceb{10}{9}{\cd}{\ca}{\cb}{\cd}
\dicea{11}{9}{\cc}{\cd}{\ca}{\cc}\diceb{11}{9}{\cc}{\ca}{\cb}{\cd}
\dicea{12}{9}{\cb}{\cb}{\cd}{\cd}\diceb{12}{9}{\cd}{\cb}{\cb}{\cc}
\dicea{13}{9}{\cd}{\ca}{\cc}{\cc}\diceb{13}{9}{\cc}{\cb}{\cb}{\cd}
\dicea{14}{9}{\cb}{\cd}{\cc}{\cd}\diceb{14}{9}{\cc}{\ca}{\ca}{\cc}
\dicea{15}{9}{\cd}{\cc}{\cc}{\cd}\diceb{15}{9}{\cc}{\ca}{\ca}{\cc}
\dicea{0}{10}{\cb}{\cb}{\ca}{\ca}\diceb{0}{10}{\cb}{\cd}{\cb}{\cc}
\dicea{1}{10}{\cd}{\cb}{\ca}{\cb}\diceb{1}{10}{\ca}{\cc}{\ca}{\cc}
\dicea{2}{10}{\ca}{\cd}{\ca}{\cb}\diceb{2}{10}{\ca}{\cc}{\ca}{\cc}
\dicea{3}{10}{\cc}{\cd}{\cb}{\cb}\diceb{3}{10}{\ca}{\cd}{\cb}{\cd}
\dicea{4}{10}{\ca}{\ca}{\cc}{\cb}\diceb{4}{10}{\cb}{\cd}{\ca}{\cd}
\dicea{5}{10}{\cc}{\cb}{\cc}{\cb}\diceb{5}{10}{\cb}{\cd}{\ca}{\cd}
\dicea{6}{10}{\ca}{\cc}{\cc}{\cb}\diceb{6}{10}{\cb}{\cd}{\ca}{\cc}
\dicea{7}{10}{\cc}{\cd}{\cd}{\ca}\diceb{7}{10}{\ca}{\cd}{\ca}{\cd}
\dicea{8}{10}{\cb}{\cb}{\ca}{\cc}\diceb{8}{10}{\cb}{\cc}{\cb}{\cc}
\dicea{9}{10}{\cd}{\ca}{\cb}{\cd}\diceb{9}{10}{\ca}{\cc}{\cb}{\cd}
\dicea{10}{10}{\cb}{\cd}{\cb}{\cd}\diceb{10}{10}{\cb}{\cc}{\cb}{\cd}
\dicea{11}{10}{\cd}{\cc}{\cb}{\cd}\diceb{11}{10}{\cb}{\cc}{\cb}{\cd}
\dicea{12}{10}{\cb}{\cb}{\cc}{\cd}\diceb{12}{10}{\cb}{\cd}{\cb}{\cc}
\dicea{13}{10}{\cd}{\cb}{\cd}{\cd}\diceb{13}{10}{\ca}{\cc}{\ca}{\cc}
\dicea{14}{10}{\ca}{\cd}{\cd}{\cd}\diceb{14}{10}{\ca}{\cc}{\ca}{\cc}
\dicea{15}{10}{\cc}{\cd}{\cd}{\cc}\diceb{15}{10}{\ca}{\cd}{\cb}{\cd}
\dicea{0}{11}{\cb}{\cb}{\cb}{\cb}\diceb{0}{11}{\cd}{\cd}{\cb}{\cc}
\dicea{1}{11}{\cd}{\ca}{\ca}{\ca}\diceb{1}{11}{\cc}{\cd}{\cb}{\cd}
\dicea{2}{11}{\cb}{\cd}{\ca}{\cb}\diceb{2}{11}{\cc}{\cc}{\ca}{\cc}
\dicea{3}{11}{\cd}{\cc}{\ca}{\cb}\diceb{3}{11}{\cc}{\cc}{\ca}{\cc}
\dicea{4}{11}{\cb}{\cb}{\cd}{\ca}\diceb{4}{11}{\cd}{\cd}{\ca}{\cc}
\dicea{5}{11}{\cd}{\cb}{\cc}{\cb}\diceb{5}{11}{\cd}{\cd}{\ca}{\cd}
\dicea{6}{11}{\ca}{\cd}{\cc}{\cb}\diceb{6}{11}{\cd}{\cd}{\ca}{\cd}
\dicea{7}{11}{\cc}{\cd}{\cc}{\cb}\diceb{7}{11}{\cc}{\cd}{\ca}{\cd}
\dicea{8}{11}{\cb}{\ca}{\cb}{\cd}\diceb{8}{11}{\cd}{\cc}{\cb}{\cd}
\dicea{9}{11}{\cd}{\ca}{\ca}{\cc}\diceb{9}{11}{\cc}{\cc}{\cb}{\cd}
\dicea{10}{11}{\ca}{\cc}{\cb}{\cd}\diceb{10}{11}{\cd}{\cc}{\cb}{\cc}
\dicea{11}{11}{\cc}{\cc}{\cb}{\cd}\diceb{11}{11}{\cd}{\cc}{\cb}{\cd}
\dicea{12}{11}{\cb}{\cb}{\cd}{\cc}\diceb{12}{11}{\cd}{\cd}{\cb}{\cc}
\dicea{13}{11}{\cd}{\ca}{\cc}{\cd}\diceb{13}{11}{\cc}{\cd}{\cb}{\cd}
\dicea{14}{11}{\cb}{\cd}{\cd}{\cd}\diceb{14}{11}{\cc}{\cc}{\ca}{\cc}
\dicea{15}{11}{\cd}{\cc}{\cd}{\cd}\diceb{15}{11}{\cc}{\cc}{\ca}{\cc}
\dicea{0}{12}{\cb}{\cb}{\cb}{\cb}\diceb{0}{12}{\cb}{\cb}{\cd}{\cc}
\dicea{1}{12}{\cd}{\cb}{\cb}{\ca}\diceb{1}{12}{\ca}{\ca}{\cc}{\cc}
\dicea{2}{12}{\ca}{\cd}{\cb}{\ca}\diceb{2}{12}{\ca}{\ca}{\cc}{\cc}
\dicea{3}{12}{\cc}{\cd}{\ca}{\ca}\diceb{3}{12}{\ca}{\cb}{\cd}{\cd}
\dicea{4}{12}{\ca}{\ca}{\cd}{\ca}\diceb{4}{12}{\cb}{\cb}{\cc}{\cd}
\dicea{5}{12}{\cc}{\cb}{\cd}{\ca}\diceb{5}{12}{\cb}{\cb}{\cc}{\cd}
\dicea{6}{12}{\ca}{\cc}{\cd}{\ca}\diceb{6}{12}{\cb}{\cb}{\cc}{\cc}
\dicea{7}{12}{\cc}{\cd}{\cc}{\cb}\diceb{7}{12}{\ca}{\cb}{\cc}{\cd}
\dicea{8}{12}{\cb}{\cb}{\cb}{\cd}\diceb{8}{12}{\cb}{\ca}{\cd}{\cc}
\dicea{9}{12}{\cd}{\ca}{\ca}{\cc}\diceb{9}{12}{\ca}{\ca}{\cd}{\cd}
\dicea{10}{12}{\cb}{\cd}{\ca}{\cc}\diceb{10}{12}{\cb}{\ca}{\cd}{\cd}
\dicea{11}{12}{\cd}{\cc}{\ca}{\cc}\diceb{11}{12}{\cb}{\ca}{\cd}{\cd}
\dicea{12}{12}{\cb}{\cb}{\cd}{\cc}\diceb{12}{12}{\cb}{\cb}{\cd}{\cc}
\dicea{13}{12}{\cd}{\cb}{\cc}{\cc}\diceb{13}{12}{\ca}{\ca}{\cc}{\cc}
\dicea{14}{12}{\ca}{\cd}{\cc}{\cc}\diceb{14}{12}{\ca}{\ca}{\cc}{\cc}
\dicea{15}{12}{\cc}{\cd}{\cc}{\cd}\diceb{15}{12}{\ca}{\cb}{\cd}{\cd}
\dicea{0}{13}{\cb}{\cb}{\ca}{\ca}\diceb{0}{13}{\cd}{\cb}{\cd}{\cc}
\dicea{1}{13}{\cd}{\ca}{\cb}{\cb}\diceb{1}{13}{\cc}{\cb}{\cd}{\cd}
\dicea{2}{13}{\cb}{\cd}{\cb}{\ca}\diceb{2}{13}{\cc}{\ca}{\cc}{\cc}
\dicea{3}{13}{\cd}{\cc}{\cb}{\ca}\diceb{3}{13}{\cc}{\ca}{\cc}{\cc}
\dicea{4}{13}{\cb}{\cb}{\cc}{\cb}\diceb{4}{13}{\cd}{\cb}{\cc}{\cc}
\dicea{5}{13}{\cd}{\cb}{\cd}{\ca}\diceb{5}{13}{\cd}{\cb}{\cc}{\cd}
\dicea{6}{13}{\ca}{\cd}{\cd}{\ca}\diceb{6}{13}{\cd}{\cb}{\cc}{\cd}
\dicea{7}{13}{\cc}{\cd}{\cd}{\ca}\diceb{7}{13}{\cc}{\cb}{\cc}{\cd}
\dicea{8}{13}{\cb}{\ca}{\ca}{\cc}\diceb{8}{13}{\cd}{\ca}{\cd}{\cd}
\dicea{9}{13}{\cd}{\ca}{\cb}{\cd}\diceb{9}{13}{\cc}{\ca}{\cd}{\cd}
\dicea{10}{13}{\ca}{\cc}{\ca}{\cc}\diceb{10}{13}{\cd}{\ca}{\cd}{\cc}
\dicea{11}{13}{\cc}{\cc}{\ca}{\cc}\diceb{11}{13}{\cd}{\ca}{\cd}{\cd}
\dicea{12}{13}{\cb}{\cb}{\cc}{\cd}\diceb{12}{13}{\cd}{\cb}{\cd}{\cc}
\dicea{13}{13}{\cd}{\ca}{\cd}{\cc}\diceb{13}{13}{\cc}{\cb}{\cd}{\cd}
\dicea{14}{13}{\cb}{\cd}{\cc}{\cc}\diceb{14}{13}{\cc}{\ca}{\cc}{\cc}
\dicea{15}{13}{\cd}{\cc}{\cc}{\cc}\diceb{15}{13}{\cc}{\ca}{\cc}{\cc}
\dicea{0}{14}{\cb}{\cb}{\ca}{\cb}\diceb{0}{14}{\cb}{\cd}{\cd}{\cc}
\dicea{1}{14}{\cd}{\cb}{\cb}{\cb}\diceb{1}{14}{\ca}{\cc}{\cc}{\cc}
\dicea{2}{14}{\ca}{\cd}{\cb}{\cb}\diceb{2}{14}{\ca}{\cc}{\cc}{\cc}
\dicea{3}{14}{\cc}{\cd}{\cb}{\ca}\diceb{3}{14}{\ca}{\cd}{\cd}{\cd}
\dicea{4}{14}{\cb}{\cb}{\cc}{\cb}\diceb{4}{14}{\cb}{\cd}{\cc}{\cc}
\dicea{5}{14}{\cd}{\ca}{\cd}{\ca}\diceb{5}{14}{\ca}{\cd}{\cc}{\cd}
\dicea{6}{14}{\cb}{\cd}{\cd}{\ca}\diceb{6}{14}{\cb}{\cd}{\cc}{\cd}
\dicea{7}{14}{\cd}{\cc}{\cd}{\ca}\diceb{7}{14}{\cb}{\cd}{\cc}{\cd}
\dicea{8}{14}{\ca}{\ca}{\cb}{\cd}\diceb{8}{14}{\cb}{\cc}{\cd}{\cd}
\dicea{9}{14}{\cc}{\cb}{\cb}{\cd}\diceb{9}{14}{\cb}{\cc}{\cd}{\cd}
\dicea{10}{14}{\ca}{\cc}{\cb}{\cd}\diceb{10}{14}{\cb}{\cc}{\cd}{\cc}
\dicea{11}{14}{\cc}{\cd}{\ca}{\cc}\diceb{11}{14}{\ca}{\cc}{\cd}{\cd}
\dicea{12}{14}{\cb}{\cb}{\cd}{\cd}\diceb{12}{14}{\cb}{\cd}{\cd}{\cc}
\dicea{13}{14}{\cd}{\cb}{\cd}{\cc}\diceb{13}{14}{\ca}{\cc}{\cc}{\cc}
\dicea{14}{14}{\ca}{\cd}{\cd}{\cc}\diceb{14}{14}{\ca}{\cc}{\cc}{\cc}
\dicea{15}{14}{\cc}{\cd}{\cc}{\cc}\diceb{15}{14}{\ca}{\cd}{\cd}{\cd}
\dicea{0}{15}{\cb}{\cb}{\cb}{\ca}\diceb{0}{15}{\cd}{\cd}{\cd}{\cc}
\dicea{1}{15}{\cd}{\ca}{\ca}{\cb}\diceb{1}{15}{\cc}{\cd}{\cd}{\cd}
\dicea{2}{15}{\cb}{\cd}{\cb}{\cb}\diceb{2}{15}{\cc}{\cc}{\cc}{\cc}
\dicea{3}{15}{\cd}{\cc}{\cb}{\cb}\diceb{3}{15}{\cc}{\cc}{\cc}{\cc}
\dicea{4}{15}{\cb}{\ca}{\cd}{\ca}\diceb{4}{15}{\cd}{\cd}{\cc}{\cd}
\dicea{5}{15}{\cd}{\ca}{\cc}{\cb}\diceb{5}{15}{\cc}{\cd}{\cc}{\cd}
\dicea{6}{15}{\ca}{\cc}{\cd}{\ca}\diceb{6}{15}{\cd}{\cd}{\cc}{\cc}
\dicea{7}{15}{\cc}{\cc}{\cd}{\ca}\diceb{7}{15}{\cd}{\cd}{\cc}{\cd}
\dicea{8}{15}{\cb}{\cb}{\ca}{\cc}\diceb{8}{15}{\cd}{\cc}{\cd}{\cc}
\dicea{9}{15}{\cd}{\cb}{\cb}{\cd}\diceb{9}{15}{\cd}{\cc}{\cd}{\cd}
\dicea{10}{15}{\ca}{\cd}{\cb}{\cd}\diceb{10}{15}{\cd}{\cc}{\cd}{\cd}
\dicea{11}{15}{\cc}{\cd}{\cb}{\cd}\diceb{11}{15}{\cc}{\cc}{\cd}{\cd}
\dicea{12}{15}{\cb}{\cb}{\cc}{\cc}\diceb{12}{15}{\cd}{\cd}{\cd}{\cc}
\dicea{13}{15}{\cd}{\ca}{\cd}{\cd}\diceb{13}{15}{\cc}{\cd}{\cd}{\cd}
\dicea{14}{15}{\cb}{\cd}{\cd}{\cc}\diceb{14}{15}{\cc}{\cc}{\cc}{\cc}
\dicea{15}{15}{\cd}{\cc}{\cd}{\cc}\diceb{15}{15}{\cc}{\cc}{\cc}{\cc}

\end{tikzpicture}
\caption{Illustration of a clique of 256 vertices in $G_{8,2}$. Each ``dice'' with eight dots represents a vertex, and each dot represents a
coordinate. A 
\ifCOLOR
black, dark blue, white, and light blue dot 
\fi
represent a coordinate set to 0, 1, 2, and 3, respectively.}
\label{fig:clique256}
\end{figure}

\section{Conclusions and Future Work}\label{sec:conclusions}

In this paper, we analyzed maximal cliques in the graphs $G_{7,3}$, $G_{7,4}$, and $G_{7,6}$ by combining symmetry-breaking and SAT-solving techniques. For the initial symmetry breaking we adapt some of the arguments of Perron. Additional symmetry breaking is performed on the propositional level and this part is mechanically verified. We partitioned the resulting formulas into thousands of subformulas and used a SAT solver to check that each subformula cannot be extended to a clique of size $128$. Additionally, we optimized and certified the resulting proofs of unsatisfiability. As a result, we proved Theorem~\ref{thm:main}, which resolves Keller's conjecture in dimension 7.

In the future, we hope to construct a formally-verified argument for 
Keller's conjecture, starting with a formalization of Keller's conjecture \ifIJCAR\else \ref{conj:keller} \fi down to the relation of the existence of cliques of size $2^n$ in Keller graphs and finally the correctness of the presented encoding. This effort would likely involve formally verifying most of the theory discussed in the Appendix\ifIJCAR of the extended version of the paper\fi\change{.} On top of that, we would like to construct a single proof of unsatisfiability that incorporates all the clausal symmetry breaking and the proof of unsatisfiability of all the subformulas and validate this proof using a formally-verified checker.

Furthermore, we would like to extend the analysis to $G_{7,s}$, including computing the size of the largest cliques for various values of $s$. Another direction to consider is to study the maximal cliques in $G_{8,s}$ in order to have some sort of classification of all maximal cliques.

\ifIJCAR
\else

\appendix

\section{Appendix}\label{app}
\normalsize

In this Appendix, we formalize the connection between Keller graphs and Keller's original conjecture on cube tilings.  In particular, we give an overview of the various results since the mid-twentieth century which have contributed toward understanding Keller's conjecture. For similar discussions and surveys on Keller's conjecture, see for example~\cite{shorminkowski,szabo1993cube}. In addition, the survey of Zong~\cite{zong2005known} situates Keller's conjecture in the broader context of research on unit cubes. An additional motivation for this appendix is to give insight to the steps needed in order to fully, formally verify the proof of Keller's conjecture in dimension 7. \changeprime{Since the first version of this appendix appeared, significant progress has been made toward formalization in Lean~\cite{clune2023formalized}.}

\subsection{Cube Tilings}

Let $d \ge 1$ be the dimension--the case $d=7$ is of the most interest in this paper. A \emph{unit cube} (or just \emph{cube}) is a translation of $[0, 1)^d$. In particular, for any $x \in \mathbb R^d$, we define $[0, 1)^d + x$ to be the translated cube $[x_1, x_1+1) \times \cdots [x_d, x_d+1)$. We call $x$ the \emph{corner} of the cube $[0, 1)^d + x$. For brevity, we denote $C^d(x) := [0, 1)^d + x$. We say that two cubes are \emph{disjoint} if they do not intersect. Disjointness is equivalent to the cube corners being ``far apart'' in some coordinate.

\propdisj*

\begin{proof}
If for some $i$, $|x_i - y_i| \ge 1$, then $[x_i, x_i+1)$ and $[y_i, y_i+1)$ cannot intersect. Therefore, $C^d(x)$ and $C^d(y)$ are disjoint.

We prove the contrapositive to show the reverse implication. Assume that for all $i$, $|x_i - y_i| < 1$. In this case, $[x_i, x_i+1)$ and $[y_i, y_i+1)$ intersect non-trivially for all $i$. Therefore, $C^d(x)$ and $C^d(y)$ intersect non-trivially.
\hfill $\square$
\end{proof}

We say that disjoint cubes $C^d(x)$ and $C^d(y)$ are \emph{facesharing} if there exists exactly one coordinate $i \in [d]$ such that $|x_i - y_i| = 1$ and for all $j \in [d]$ such that $j \neq i$, $x_j = y_j$. This is equivalent to saying that $x = y \pm e_i$, where $e_1, \hdots, e_d$ are the standard unit basis vectors.

Let $T \subset \mathbb R^d$ be a set of cube corners. $T$ is a \emph{cube tiling} of $\mathbb{R}^d$ if $[0,1)^d + T = \{C^d(t) : t \in T\}$ is a family of pairwise disjoint cubes such that $\bigcup_{t\in T} C^d(t) = \mathbb{R}^d$. Note that $T = \mathbb Z^d$ produces the standard lattice 
cube tiling of $\mathbb R^d$.

We say that a cube tiling is \emph{faceshare-free} if no pair of distinct cubes in the tiling share a face. For example, $[0, 1)^d + \mathbb Z^d$ is not faceshare-free since, e.g., $[0, 1)^d$ and $[0, 1)^d+e_1$ faceshare. Keller conjectured that all tilings have a facesharing pair of cubes. %

\kellerconj*

Perron~\cite{perron1940luckenlose,perron1940luckenlose2} was able to resolve the conjecture in dimensions $d \le 6$ using combinatorial casework\change{. Lagarias and Shor~\cite{lagarias1992keller} and Mackey~\cite{mackey2002cube} showed the conjecture to be false for dimensions $d \ge 8$. This paper resolves the final case of $d=7$ in the affirmative.}

\change{We shall explain in the remainder of this Appendix the key structural results in the literature that led to the full resolution of Keller's conjecture.}

\subsection{Structure of Tilings}\label{subsec:structure}

Let $T \subset \mathbb R^d$ be any set of cube corners such that $[0, 1)^d + T$ is a cube tiling. By the definition of cube tiling, we know that for any $x \in \mathbb Z^d$, there exists a unique $t \in T$ such that $x \in C^d(t)$. 

Likewise, for any $t \in T$, there exists a unique $x \in \mathbb Z^d$ such that $x \in C^d(t)$. This is the point $x := (\lceil t_1\rceil, \lceil t_2 \rceil, \hdots, \lceil t_n\rceil)$, where $\lceil r \rceil$ represents the least integer that $r$ does not exceed. Thus, we can refer to each cube by the integral point it contains. For all $x \in \mathbb Z^d$, we let $t(x) \in T$ be the unique corner such that $x \in C^d(t(x))$. 

For all $x \in \mathbb R^d$ and $i \in [d]$, let $\ell_i(x)$ the be line through $x$ parallel to the $i$th coordinate. Formally,
\[
\ell_i(x) = \{y \in \mathbb R^d : \forall j \in [d] \setminus \{i\}, y_j = x_j \}.
\]
We use the term \emph{line} to refer to these axis-parallel lines. We also define
\[
T_i(x) = \{t \in T : C^d(\change{t}) \cap \ell_i(x) \neq \emptyset\}
\]
We say that $T_i(x)$ is an \emph{$i$-lattice} if there exists $a \in \mathbb R$ such that for all $t \in T_i(x)$, $t_i \equiv a \mod 1$ and for all $y \in \mathbb Z$, there exists exactly one $t \in T_i(x)$ such that $t_i = a + y$. In other words, the $i$th coordinates of $T_i(x)$ are a ``shift'' of 
$\mathbb Z$.

\change{To illustrate this definition, consider the tiling $T =
  \{(a+ 0.5 b,b) + [0,1)^2 : a,b \in \mathbb Z\}.$ That is a tiling of
  cubes in $\mathbb R^2$ where the cubes line up in horizontal rows, but
  each row is shifted by $0.5$ from the next row. For all $x \in
  \mathbb R^2$, $T_1(x)$ is one of these rows of cubes. However,
  $T_2(x)$ is a ``zig-zag'' column of cubes, where each successive
  cube moves back/forth by $1/2$ across the vertical line through $x$.}

\begin{prop}\label{prop:lattice}
Let $T \subset \mathbb R^d$ be a collection of cube corners. Fix $i \in [d]$.  $[0, 1)^d + T$ is a tiling if and only if for all $x \in \mathbb R^d$, $T_i(x)$ is an $i$-lattice.
\end{prop}

\begin{proof}
First, assume $[0, 1)^d + T$ is a tiling. Note that the intersection of any unit cube with any line is either the empty set or a half-open interval of length $1$. Furthermore, the half-open intervals must partition the line, so the starting points of the intervals are integrally spaced. Thus, $T_i(x)$ is an $i$-lattice.

Likewise, if $T_i(x)$ is an $i$-lattice for all $x \in \mathbb R^d$, then every point on each line is in exactly one cube. Thus, every point is in exactly one cube, so $[0, 1)^d + T$ is a tiling.
\hfill $\square$
\end{proof}

The following structural result also follows.

\begin{prop}\label{prop:buddy}
Let $T \subset \mathbb R^d$ be such that $[0, 1)^d + T$ is a cube tiling. For all $x \in \mathbb Z^d$ and $i \in [d]$,
\[
t(x)_i + 1 = t(x + e_i)_i.
\]
\end{prop}

\begin{proof}
Consider the line $\ell_i(x)$. The intervals $C^d(t(x)) \cap \ell_i(x)$ and $C^d(t(x+e_i)) \cap \ell_i(x)$ must be adjacent since they contain adjacent integral points. Therefore, by Proposition~\ref{prop:lattice}, $t(x)_i + 1 = t(x+e_i)_i$.
\hfill $\square$
\end{proof}

The $t(x)$ notation also gives a straightforward condition for checking whether the corresponding tiling is faceshare-free.

\begin{prop}\label{prop:faceshare}
Let $[0, 1)^d+T$ be a tiling. If for all $x\in \mathbb R^d$ and all $i \in [d]$ we have that $C^d(t(x))$ and $C^d(t(x + e_i))$ do not faceshare, then $[0, 1)^d+T$ is faceshare free.
\end{prop}

\begin{proof}
Assume for sake of contradiction some pair of cubes $C^d(t(x))$ and $C^d(t(y))$ faceshare. This implies that $t(x) = t(y) \pm e_i$ for some $i$. Since $x$ is each coordinate of $t(x)$ rounded up, and $y$ is each coordinate of $t(y)$ rounded up, $x = y \pm e_i$. This is a contradiction.
\hfill $\square$
\end{proof}

The specific values of the coordinates in a tiling are somewhat artificial. That is, the values of the coordinates can be changed while preserving the tiling, as long as the following rule is observed:

\proprepl*

\change{To partially illustrate the lemma, recall the tiling $T =
  \{(a+ 0.5 b,b) + [0,1)^2 : a,b \in \mathbb Z\}.$ By the replacement
  lemma, if we shift every other row so that it is offset by $1/3$
  instead of $1/2$,
  then no \emph{new} facesharing is introduced.}

\begin{proof}
For all $x \in \mathbb R^d$, consider the $i$-lattice $T_i(x)$. Either every $i$th coordinate in this lattice is $\not\equiv a \mod 1$, in which case nothing changes. Otherwise, every $i$th coordinate is $\equiv a \mod 1$, in which case adding $b$ to every $i$th coordinate preserves it as an $i$-lattice. Thus, by Proposition~\ref{prop:lattice}, $[0, 1)^d + T'$ is a tiling.

Now assume $[0, 1)^d + T$ is faceshare-free and there exists no $t \in T$ such that $t_i \equiv a + b\mod 1$. If $[0, 1)^d + T'$ has facesharing, there are $t^1, t^2 \in T'$ such that $t^1 = t^2 \pm e_j$ for some $j \in [d]$. In particular, $t^1_k \equiv t^2_k \mod 1$ for all $k\in[d]$.

If $t^1_i \not\equiv a+b\mod 1$, then $t^1$ and $t^2$ were not shifted when transforming $T$ to $T'$. Thus, $t^1$ and $t^2$ are a facesharing pair in $T$, \change{which is} a contradiction.

Otherwise, if $t^1_i \change{\equiv} a+b\mod 1$, then we know that $t^1$ and $t^2$ \change{were} shifted by $be_i$ when going from $T$ to $T'$ (this uses the condition that $t_i \not\equiv a + b\mod 1$ for all $t \in T$). Thus, $t^1-be_i$ and $t^2-be_i$ faceshare in $T$, again a contradiction.
\hfill $\square$
\end{proof}

\subsection{Reduction to Periodic Tilings}

In this subsection, we show \change{it} suffices to look at
\emph{periodic} tilings. We say that $T \subset \mathbb R^d$ and its
cube tiling $[0, 1)^d + T$ are \emph{periodic} if for all $t \in T$
and \change{$x \in \mathbb \mathbb Z^d$}, we also have $t + 2x \in
T$. For instance, $T = \mathbb Z^d$ is periodic, \change{as well as
  the tiling $T = \{(a+ 0.5 b,b) + [0,1)^2 : a,b \in \mathbb Z\}.$
  mentioned earlier.} In the $t(x)$ notation (see
Section~\ref{subsec:structure}), we have that for all
$x, y \in \mathbb Z^d$,
\[
t(x + 2y) = t(x) + 2y.
\]

For a given dimension $d$, if Keller's conjecture is true for all tilings then it is also true for the periodic tilings. The work of Haj\'os~\cite{hajos,hajos1950factorisation} shows that the reverse implication is also true.

\hajos*

\begin{proof}
We prove this by contradiction. Assume there exists $T \in \mathbb R^d$, not necessarily periodic, such that $[0, 1)^d + T$ is faceshare free.

Let $\hat{T} \subset T$ be the corners whose corresponding cubes \change{contain the points} $\{0, 1\}^d$. In the $t(x)$ notation, we have
\[\hat{T} := \{t(x) : x \in \{0, 1\}^d\}.\]
Note that $\hat{T}$ has exactly $2^d$ cubes (which is ultimately why the maximum clique size of the Keller graph is $2^d$).

We extend $\hat{T}$ to a new periodic set $T'$ as follows
\[
    T' := \{t' = t + 2x \mid t \in \hat{T}, x \in \mathbb Z^d\}.
\]

Observe that each point $y \in \mathbb Z^d$ is in exactly one cube of $[0, 1)^d + T'$. Call this cube $t'(y)$.

We claim that $[0, 1)^d + T'$ is a periodic faceshare-free tiling of 
$\mathbb R^d$. This follows from the following three facts:

\begin{enumerate}
\item \emph{Every point of $[0, 1]^d$ is in exactly one cube of $[0, 1)^d + T'$.} Since $[0, 1)^d + T$ is a tiling, every point of $[0, 1]^d$ is in exactly one cube of $T$. Furthermore, any cube which covers a point of $[0, 1]^d$ must cover one of the corners. Since the cubes which cover $\{0, 1\}^d$ are the same for $T$ and $T'$ (by definition), every point of $[0, 1]^d$ must also be in exactly one cube of $[0, 1)^d + T'$.

\item \emph{$[0, 1)^d + T'$ is a tiling.} We seek to show that for all $i \in [d]$, every point of the closed cube $[0, 1]^d + e_i$ is in exactly one point of the tiling $[0, 1)^d + T'$. Fix $x \in [0, 1]^d + e_i$. Let $x^{i\to 0}$ be the point identical to $x$ except the $i$th coordinate is replaced with $0$. Define $x^{i\to 1}$ and $x^{i\to 2}$ similarly. Let $\ell_1$ be the closed line segment from $x^{i\to 0}$ to $x^{i\to 1}$ and $\ell_2$ be the closed line segment from $x^{i\to 1}$ to $x^{i\to 2}$. 

Since every point of $[0, 1]^d$ is in exactly one cube of $T'$, there must exist $y, z \in \{0, 1\}^d$ such that $x^{i\to 0} \in C^d(t'(y))$ and $x^{i\to 1} \in C^d(t'(z))$. By an argument similar to that of Proposition~\ref{prop:buddy}, this implies that $\ell_1$ is covered by the union of $C^d(t'(y))$ and $C^d(t'(z))$ and that $t'(y)_i+1 = t'(z)_i$. 

Now, consider $y' = y + 2e_i$, observe that since $x^{i\to 0} \in C^d(t'(y))$, $x^{i\to 0} +2e_2 = x^{i\to 2} \in C^d(t'(y'))$. Also, $t'(z)_i + 1 = t'(y')_i$. Therefore, $C^d(t'(z))$ and $C^d(t'(y'))$ disjointly cover $\ell_2$ and thus exactly one of them contains $x$.

To see that no other cube can cover $x$, assume that a third cube $C^d(t'(w))$ covers $\ell_2$. Then, either a portion of $C^d(t'(w))$ intersects $\ell_1$ or a portion of $C^d(t'(w-2e_i))$ intersects $\ell_1$. In either case, since every point of $\ell_1$ is covered by exactly one cube, we deduce that $w = z$ or $w = y + 2e_i$. Thus, $x$ is covered by exactly one cube.

By a similar argument, every point of $[0, 1]^d - e_i$ is covered by exactly one point of $[0, 1)^d + T'$. By a suitable induction, every closed cube $[0, 1]^d + y$ for some $y \in \mathbb Z^d$ has every point covered by exactly one cube of $[0, 1)^d + T'$. Thus, every point of $\mathbb R^d$ is covered by exactly one cube of $[0, 1)^d + T'$, and thus $[0, 1)^d + T'$ is a tiling.

\item \emph{$[0, 1)^d + T'$ is faceshare free.} If not, then by Proposition~\ref{prop:faceshare}, we must have that $t'(x) +e_i= t'(x+e_i)$ for some $x \in \mathbb Z^d$ and $i \in [d]$. Note there exist $x' \in \{0, 1\}^d$ and $y \in \mathbb Z^d$ such that $x = x' + 2y$. Then,
\begin{align*}
t'(x') + e_i &= t'(x'+2y)-2y + e_i\\
&= t'(x)+e_i - 2y\\
&= t'(x+e_i) - 2y\\
&= t'(x'+e_i).
\end{align*}
If $x'_i = 0$, then $x'+e_i \in \{0, 1\}^n$ as well. Thus, $C^d(t'(x')), C^d(t'(x'+e_i))$ faceshare in $[0,1)^d+\hat{T}$ (and thus $[0, 1)^d+T$), \change{which is} a contradiction.

If $x'_i=1$, then, $x'-e_i \in \{0, 1\}^n$. And a similar calculation shows that $t'(x' - e_i) + e_i = t'(x')$. Therefore, $C^d(t'(x')), C^d(t'(x'-e_i))$ faceshare in $\hat{T}$ (and $T$), which is a contradiction.
\hfill $\square$
\end{enumerate}
\end{proof}

\subsection{Reduction to Keller graphs}

In the previous section, we showed that Keller's conjecture for $\mathbb R^d$ is equivalent to showing that there is no faceshare-free periodic tiling of $\mathbb R^d$. Now we show faceshare-free periodic tilings correspond to cliques in the Keller graph.

We say that a periodic tiling is \emph{$s$-discrete} if every coordinate has at most $s$ distinct values modulo $1$ in the tiling. A key observation is that $s$ is bounded by the a function of the dimension.\footnote{This is a (ultimately equivalent) variant of a definition of Kisielewicz and {\L}ysakowska~\cite{kisielewicz2017rigid,kisielewicz2015keller,kisielewicz2017towards}. They defined for $x \in  \mathbb{R}^d$ and $i \in [d]$, $L(T,x,i)$ to be the set of all $i^{\text{th}}$ coordinates $t_i$ of vectors $t \in T$ such that 
$([0,1)^d + t) \cap ([0,1]^d + x) \neq \emptyset$ and $t_i \leq x_i$. A tiling is $s$-discrete if and only if $L(T, x, i) \le s$ for all $x$ and $i$.}

\propszabo*

\begin{proof}
For a given periodic tiling $T$, every cube is an integral translate of one of the $2^d$ cubes $t(x)$ for $x \in \{0, 1\}^d$. By Proposition~\ref{prop:buddy}, we know that for all $i \in [d]$, $t(x)_i$ and $t(x \oplus e_i)_i$ are the same modulo $1$. Pairing up the elements of $\{0, 1\}^d$, we see that there can be at most $2^{d-1}$ distinct values modulo 1 in the $i$th coordinate.
\hfill $\square$
\end{proof}

\begin{lem}
$G_{d,s}$ has a clique of size $2^d$ if and only if there exists a faceshare-free $s$-discrete periodic tiling $[0, 1)^d+T$ in dimension $d$.
\end{lem}

\begin{proof}
First, let $[0, 1)^d + T$ be an $s$-discrete periodic tiling. By repeatedly applying the Replacement Lemma, we can replace each coordinate of the entries of $T$ with elements of $\left\{0, \dfrac{1}{s}, \hdots, \dfrac{s-1}{s}\right\}$ modulo $1$. In other words, we may assume without loss of generality that $T \subset (\frac{1}{s}\mathbb Z)^d$.

Note that each $t(x) \in T$ $(x \in \mathbb Z^d)$ corresponds to a vertex in $G_{d, s}$. To see why, note that $st(x) \in \mathbb Z^d$. Thus, there is a unique $u(x) \in \{0, 1, \hdots, 2s-1\}^d$ (the vertices of $G_{d, s}$) such that $st(x) \equiv u(x) \mod 2s$.

We claim that $K = \{u(x) : x \in \{0, 1\}^d\}$ is a clique of $G_{d,s}$.  Assume for sake of contradiction that for some $x, y \in \{0, 1\}^d$, $u(x)$ and $u(y)$ are not connected by an edge. We have two cases to consider
\begin{enumerate}
\item \emph{$u(x)_i \neq u(y)_i+s$ for all $i \in [d]$.} Thus, $t(x)_i \not\equiv t(y)_i + 1 \mod 2$ for all $i \in [d]$. Then for all $i \in [d]$, there exists a positive integer $z_i \in \mathbb Z$ such that $[t(x)_i + 2z, t(x_i) + 2z+1)$ intersects $[t(y)_i, t(y)_i+1)$. In other words, $t(x+2z)$ and $t(y)$ intersect non-trivially. This is a contradiction.

\item \emph{$u(x)_i = u(y)_i+s$ for some $i \in [d]$ and $u(x)_j = u(y)_j$ for all $j \neq i$.} Then, $t(x) = t(y)+e_i$. This contradicts that $T$ is faceshare-free.
\end{enumerate}
Thus, $K$ is indeed a clique of size $2^d$ in $G_{n, s}$ (note that case 1 above rules out repeated vertices).

For the reverse implication, let $K \subset \{0, 1, \hdots, 2s-1\}^d$ be a clique of size $2^d$. Define a periodic set of cube corners $T$ as follows
\[
    T := \left\{\frac{u}{s} + 2x : u \in K, x \in \mathbb Z^d\right\}.
\]
Note that $T$ is $s$-discrete because modulo $1$ each coordinate must be in the set $\{0, \frac{1}{s}, \hdots, \frac{s-1}{s}\}$. 

Next we show every pair of cubes in $[0, 1)^d + T$ is disjoint. If not, by the definition of $T$, there exist $u, v \in K$ and $x, y \in \mathbb Z^d$ such that $C^d(u/s + 2x)$ and $C^d(v/s+ 2y)$ intersect. By Proposition~\ref{prop:disj}, this means for all $i \in [d]$, 
\[
\left|\frac{u_i - v_i}{s} + 2x - 2y\right| < 1.
\]
If $u = v$, then we must have $x = y$, so the cubes are identical. Otherwise, the above inequality implies that $u_i \neq v_i \pm s$ for all $i \in [d]$, which means that $u$ and $v$ cannot be in the same clique, \change{which is} a contradiction.

To show that $[0, 1)^d + T$ is a tiling, it suffices to show that every point $x \in \mathbb R^d$ is in some cube of the tiling. We partition $\mathbb R^d$ into \emph{cubelets}: $[0, 1/s)^d + (\frac{1}{s}\mathbb R)^d$. It suffices to show that each cubelet is in some cube of the tiling. In fact, each cube of $T$ contains $s^d$ cubelets. Call two cubelets \emph{equivalent} if their corners differ by an even integer vector. For each coordinate, there are $2s$ non-equivalent starting coordinates $\{0, \frac{1}{s}, \hdots, \frac{2s-1}{s}\}$ so there are a total of $(2s)^d$ inequivalent cubelets.

Observe that for all $u \in K$ and $x, y \in \mathbb Z^d$, $C^d(\frac{u}{s} + x)$ and $C^d(\frac{u}{s}+2y)$ cover equivalent cubelets. Another observation is that if $u, v \in K$ are distinct then $C^d(\frac{u}{s})$ and $C^d(\frac{v}{s})$ cover inequivalent cubelets. Otherwise, by the equivalence we could shift one of the cubes by an even integer vector to make them overlap, which contradicts the assumption that the cubes are disjoint. 
Since $K$ has size $2^d$, this means $[0, 1)^d + T$ covers at least 
$2^ds^d = (2s)^d$ inequivalent cubelets. This means that all cubelets are covered; therefore all of $\mathbb R^d$ is covered.

Next, we show that $[0, 1)^d + T$ is faceshare-free. Otherwise, there exists, $u, v\in K$ and $x, y \in \mathbb Z^d$ and $i \in [d]$ such that
\[
\frac{u}{s} + 2x = \frac{v}{s} + 2y + e_i.
\]
This would imply that $u_i = v_i \pm s$ and $u_j = v_j$ for all $j \neq i$, \change{which is} a contradiction.
\hfill $\square$
\end{proof}

As a corollary, we deduce that Keller's conjecture is equivalent to the maxclique problem for a suitable Keller graph.

\corradiszabo*

The relatively recent line of papers by Kisielewicz and {\L}ysakowska improved on the above theorem in dimension $7$ by showing that any potential faceshare-free tilings must be $s$-discrete for a small value of $s$ (in comparison to $64$). These results are stated as follows.

\begin{thm}[\cite{kisielewicz2017rigid}]
Every periodic faceshare-free tiling in $7$ dimensions is $6$-discrete. Therefore, Keller's conjecture in dimension $7$ is equivalent to the lack of a clique of size $128$ in $G_{7,6}$.
\end{thm}

\begin{thm}[\cite{kisielewicz2015keller}]
Every periodic faceshare-free tiling in $7$ dimensions is $4$-discrete. Therefore, Keller's conjecture in dimension $7$ is equivalent to the lack of a clique of size $128$ in $G_{7,4}$.
\end{thm}

\begin{thm}[Corollary 1.3 of~\cite{kisielewicz2017towards}]
Every periodic faceshare-free tiling in $7$ dimensions is $3$-discrete. Therefore, Keller's conjecture in dimension $7$ is equivalent to the lack of a clique of size $128$ in $G_{7,3}$.
\end{thm}
\fi

\section*{Acknowledgments}
The authors acknowledge the Texas Advanced Computing Center
(TACC) %
at The University of Texas at Austin, RIT Research Computing, and the
Stanford Research Computing Center for providing HPC resources that
have contributed to the research results reported within this
paper. Joshua is supported by an NSF graduate research
fellowship. Marijn and David are supported by NSF grant
CCF-2006363. David is supported in part by NSF grant CCF-2030859 to
the Computing Research Association for the CIFellows Project. We thank
Andrzej Kisielewicz and Jasmin Blanchette for valuable comments on an
earlier version of the manuscript. We thank William Cooperman for
helpful discussions on a previous attempt at programming simulations
to study the half-integral case. We thank Alex Ozdemir for helpful
feedback on both the paper and the codebase. We thank Xinyu Wu for
making this collaboration possible. \change{We thank Joshua Clune for
  pointing out an error in the manuscript, and anonymous reviewers for
  much helpful feedback.}

\bibliographystyle{spbasic}
\bibliography{ref}

\end{document}